\documentclass[12pt]{amsart}

\usepackage{amsmath}
\usepackage{amssymb}

\newtheorem{theorem}{Theorem}[section]

\newtheorem*{main1}{Theorem I}

\newtheorem{lemma}[theorem]{Lemma}
\newtheorem{proposition}[theorem]{Proposition}
\newtheorem{corollary}[theorem]{Corollary}

\theoremstyle{definition}

\theoremstyle{remark}
\newtheorem{remark}{Remark}

\newcommand{\Proof}{\begin{proof}} 

\def\dots{\mathinner{\ldotp\ldotp\ldotp}}

\def\R{\mathbb R}

\def\C{\mathbb C}
\def\D{\mathbf D}

\numberwithin{equation}{section}
\allowdisplaybreaks
\begin{document}

\title[Form boundedness of the second order differential operator]
{Form boundedness of  the general\\ second order  
differential operator}

\author[V.~G. Maz'ya]
{V.~G. Maz'ya}
\address{Department of Mathematics,
Ohio State University, 231 West 18-th
Ave.,  Columbus, OH 43210, USA}
\address{Department of Mathematics,
Link\"oping University, SE-581 83, Link\"oping, 
Sweden}
\email{vlmaz@mai.liu.se}

\author[I.~E. Verbitsky]
{I.~E. Verbitsky}
\address{Department of Mathematics,
University of Missouri,
Columbia, MO 65211, USA}
\email{igor@math.missouri.edu}

\begin{abstract} 
 
 We give explicit  necessary and sufficient conditions 
for the   boundedness of the  general 
 second order differential operator 
$$\mathcal L = \sum_{i, \, j=1}^n \, a_{ij} \, \partial_i \partial_j + 
\sum_{j=1}^n \, b_{j} \,  \partial_j + c$$
 with real- or  complex-valued 
distributional coefficients $a_{ij}$, $b_{j}$, and $c$, acting from the Sobolev 
space $W^{1, \, 2}(\R^n)$ 
 to its dual $W^{-1, \, 2}(\R^n)$.   This enables us 
to obtain analytic criteria for the fundamental notions of relative form boundedness, compactness, and 
infinitesimal form boundedness of $\mathcal{L}$ with respect to the Laplacian on $L^2(\R^n)$. 

In particular, we establish a complete characterization 
of  the  form boundedness  of 
 the Schr\"odinger operator  $(i \nabla \, + \vec a)^2 + q$  
with magnetic vector potential $\vec a \in L^2_{{\rm loc}} (\R^n)$ 
and $q \in D'(\R^n)$.

\end{abstract}

\maketitle

\section{Introduction}\label{Introduction}

The property of form boundedness, as well as related notions of  relative compactness, 
infinitesimal form 
boundedness, and subordination  of differential operators 
in Hilbert spaces, are used extensively in 
mathematical physics, geometry, and PDE, especially in relation 
to quantum mechanics problems   \cite{ChWW}, \cite{Fef}, \cite{LL}, \cite{RS},   elliptic differential 
operators and spectral theory \cite{D2},  \cite{EE},  \cite{GT}, \cite{RSS},  \cite{Sch}, \cite{Sh}, 
semigroup theory \cite{D1}, \cite{LPS}, \cite{Sim}, 
harmonic maps 
 \cite{Ev},  and Markov processes \cite{CWZ}, \cite{CrZ}. 

 The goal of the present paper is to give  an analytic 
characterization of form boundedness for the general second order 
differential operator 
\begin{equation}\label{E:1.0}
\mathcal L = \sum_{i, \, j=1}^n \, a_{ij} \, \partial_i \partial_j + 
\sum_{j=1}^n \, b_{j} \,  \partial_j + c,  
\end{equation}
where $a_{ij}$,  $b_i$, and $c$ are real- or 
complex-valued distributions, on the Sobolev space $W^{1, \, 2}(\R^n)$, and its homogeneous counterpart 
$L^{1, \, 2}(\R^n)$. 

One of our motivations is to give a  criterion for the relative form 
boundedness   of the operator $\vec b \cdot \nabla + q$ with distributional coefficients $\vec b$ and $q$  
with respect to the Laplacian  $\Delta$ on $L^2(\R^n)$. 
This ensures, in view of  the so-called 
KLMN Theorem (see \cite{EE}, Theorem IV.4.2; \cite{RS}, Theorem X.17),  
that $\mathcal L = \Delta + \vec b \cdot \nabla + q$ can be defined,  
under appropriate 
smallness assumptions on 
$\vec b$ and $q$, as an 
m-sectorial operator on $L^2(\R^n)$ so that its  quadratic form domain 
coincides with $W^{1, \, 2}(\R^n)$.  
 
 In particular, we will obtain a characterization of the relative 
form boundedness  for the magnetic Schr\"odinger operator 
\begin{equation}\label{E:1.1m}
\mathcal{M} = (i \, \nabla + \vec a )^2 + q, 
\end{equation}
with  arbitrary vector potential $\vec a \in L^2_{{\rm loc}} (\R^n)$, and 
$q \in D'(\R^n)$ on $L^2(\R^n)$ with respect to $\Delta$.

Our approach is based on factorization of functions in Sobolev spaces and  integral estimates of 
potentials of equilibrium measures, 
 combined with compensated compactness 
arguments, commutator estimates, and the idea of gauge invariance. We 
are able to treat general second order differential operators, and establish 
 an explicit Hodge decomposition for form bounded 
vector fields. It is worth mentioning that in this decomposition,  the irrotational part of the vector field 
is subject to a more stringent condition than its divergence-free counterpart. 

Methods and techniques proposed in the present paper, along with their natural extensions to 
higher order differential operators 
and more general $L^p$-inequalities,  might  be useful in further applications to mathematical physics, 
 dynamics, analysis of phases, and other nonlinear problems (see, e.g., \cite{BB1}, 
\cite{BB2}, \cite{D3}, \cite{IM}). 

For the sake of convenience,  let us assume in the Introduction 
that the principal 
part of $\mathcal L$ is  in 
the divergence form, i.e., 
\begin{equation}\label{E:1.00}
\mathcal L \, u= {\rm div} \, (A \, \nabla u) + \vec b \cdot \nabla u + q \, u,  
\qquad u \in C^\infty_0(\R^n), 
\end{equation}
where $A= (a_{ij})_{i, \, j=1}^n \in D'(\R^n)^{n \times n}$, 
$\vec b = (b_j)_{j=1}^n\in D'(\R^n)^{n}$, 
and $q \in D'(\R^n)$. 

We will present necessary and sufficient conditions on $A$, $\vec b$, 
and $q$ which guarantee the boundedness of  
 the sesquilinear form associated with $\mathcal L$:
\begin{equation}\label{E:1.1a}
| \langle \mathcal{L} \, u, \, v \rangle | \leq \, C \, ||u||_{L^{1, \, 2} (\R^n)} \, 
||v||_{L^{1, \, 2} (\R^n)}    
\end{equation}
where the constant $C$ does not depend on $u, \, v  \in C^\infty_0(\R^n)$. Here $L^{1, \, 2}(\R^n)$   is 
the completion of (complex-valued) $C^\infty_0(\R^n)$ functions with respect to the  norm 
$|| u||_{L^{1, \, 2}(\R^n)} = ||\nabla  u||_{L^2(\R^n)}$. 

Equivalently, we characterize all $A$, $\vec b$, and $q$ such that 
\begin{equation}\label{E:1.1b}
\mathcal L : \, \, L^{1, \, 2}(\R^n)\to L^{-1, \, 2}(\R^n)
\end{equation}
 is a bounded operator, where  
$L^{-1, \, 2}(\R^n)=L^{1, \, 2}(\R^n)^*$ is a dual 
Sobolev space. Analogous results are obtained below for the inhomogeneous Sobolev space 
 $W^{1, \, 2}(\R^n)= L^{1, \, 2}(\R^n)\cap L^2(\R^n)$ as well.

In the special case where  $A$, 
$\vec b$ and $q $ are locally integrable,  
the form boundedness  of $\mathcal L$ may be expressed in the form of the integral inequality 
\begin{equation}\label{E:1.1000}
\left \vert \int_{\R^n}  ( - (A \, \nabla u) \cdot \nabla \overline{v}  +  \vec b\cdot\nabla u \, \, 
\overline{v} + q  u \,  \overline{v} ) \, dx 
\right \vert  \leq C \, ||u||_{L^{1, \, 2} (\R^n)}  
||v||_{L^{1, \, 2} (\R^n)},     
\end{equation}
where the constant $C$ does not depend on  $u, \, v  \in C^\infty_0(\R^n)$.   Sometimes it will be convenient to write  
(\ref{E:1.1a}) in this form even for distributional coefficients $a_{ij}$, $b_j$, and $q $.

 To state our main results, we introduce the class of {\it admissible measures}  
 $\mathfrak{M}^{1, \, 2}_+(\R^n)$, i.e., nonnegative Borel measures $\mu$ on $\R^n$ 
 which obey  the trace inequality    
\begin{equation}\label{E:1.tr}
\int_{\R^n}  
  |u|^2  \, d \mu  \le C \, 
|| u||^2_{L^{1, \, 2}(\R^n)}, \qquad u \in C^\infty_0(\R^n), 
\end{equation} 
where the constant $C$ does not depend on $u$. For admissible measures $q (x) \, dx$ 
with nonnenegative  density $q\in L^1_{{\rm loc}}(\R^n)$, we will write $q \in \mathfrak{M}^{1, \, 2}_+(\R^n)$. 

Inequalities of this type  (with  $\mu$ possibly singular with respect to Lebesgue measure)  
have been thoroughly studied.   A straightforward consequence of (\ref{E:1.tr}) is that if   
  $\mu \in \mathfrak{M}^{1, \, 2}_+(\R^n)$ then 
\begin{equation}\label{E:1.7} 
\int_{|x-y|<r}    d \mu(y)  \le {\rm const} \, \,  r^{n-2},   
\qquad 
\end{equation}
for all $r>0, \, x \in \R^n$, if $n\ge 3$, and  $\mu =0$ if $n=1, 2$ (see e.g.  \cite{M}, Sec. 2.4). 
 
A close sufficient condition on $q\in L^1_{{\rm loc}}(\R^n)$, $q \ge 0$, 
which ensures  that $q \in \mathfrak{M}^{1, \, 2}_+(\R^n)$,   is provided by the Fefferman--Phong class 
\begin{equation}\label{E:1.5b} 
\int_{|x-y|<r}   \, q^{1 + \epsilon} 
\,   dy \le {\rm const} \, r^{n-2(1 + \epsilon) }, 
\end{equation} 
where $\epsilon > 0$, and the constant does not depend on $r>0, \, x \in \R^n$.   
More precise sufficiency results are due to Chang, Wilson, and Wolff \cite{ChWW}. 

A complete characterization of the class of admissible  measures $\mathfrak{M}^{1, \, 2}_+(\R^n)$ can be expressed in several equivalent forms:    
using capacities \cite{M}, local energy estimates \cite{KS}, pointwise 
potential inequalities \cite{MV1}, or dyadic Carleson measures 
\cite{V}. These criteria, discussed in Sec.~\ref{Section 2} below,   
employ various degrees of localization of $\mu$, and each of them has its own advantages depending on the area of application. 

We now state our main form boundedness criterion. For $A= (a_{ij})$, let  $A^t = (a_{ji})$ denote the transposed matrix, and let 
${\rm Div}\colon D'(\R^n)^{n \times n} \to D'(\R^n)$ be the row 
divergence operator  defined by 
\begin{equation}\label{E:1.Div}
{\rm Div} (a_{ij}) =  ( \sum_{j=1}^n \, \partial_j \, a_{ij} )_{i=1}^n. 
\end{equation}

\begin{main1}\label{Theorem I} Let 
$\mathcal L = {\rm div} \, (A \, \nabla\cdot) +
\vec b \cdot \nabla + q$, where  
$A \in D'(\R^n)^{n\times n}$, 
 $\vec b \in D'(\R^n)^n$ and $q \in D'(\R^n)$, $n \ge 2$.  Then 
the following statements hold. 

{\rm(i)} The sesquilinear form of $\mathcal L$ is bounded, i.e., 
 {\rm(\ref{E:1.1a})} holds 
if and only  if  $\tfrac 1 2 \, (A + A^t) \in L^\infty(\R^n)^{n \times n}$, and 
 $\, \vec b$ and $q$ can be represented respectively in the form 
\begin{equation}\label{E:1.4u}
\vec b = \vec c + {\rm Div} \, F, \qquad q = {\rm div} \, \vec h,
\end{equation} 
where $F$ 
 is a skew-symmetric matrix field such that
\begin{equation}\label{E:1.4v}
F - \tfrac 12 \, (A-A^t) \in {\rm BMO(\R^n)}^{n\times n},
\end{equation} 
whereas $\vec c$ and $\vec h$ belong to $L^2_{{\rm loc}}(\R^n)^{n}$, and 
obey the condition 
\begin{equation}\label{E:1.5u}
|\vec c|^2 + | \vec h |^2  \in \mathfrak{M}^{1, \, 2}_+(\R^n). 
\end{equation}

{\rm(ii)} If  the sesquilinear form of $\mathcal L$ is bounded, then $\vec c$,  
$F$, and $\vec h$  in decomposition {\rm (\ref{E:1.4u})} can be determined  
explicitly by    
\begin{align}\label{E:1.6u}
\vec c & = \nabla ( \Delta^{-1} {\rm div} \, \vec b), 
 \qquad \vec h = 
 \nabla (\Delta^{-1} \, q), \\
\qquad F & = \Delta^{-1} {\rm curl} \, [ \vec b - \tfrac 1 2 \,  {\rm Div} \, (A-A^t)] + \tfrac 1 2 (A-A^t). 
\label{E:1.6uu}
\end{align} 
where 
\begin{equation}\label{E:1.6v}
\Delta^{-1} {\rm curl} \, [ \vec b - \tfrac 1 2 \, 
{\rm Div} \, (A-A^t)]\in {\rm BMO}(\R^n)^{n \times n}, 
\end{equation} 
and 
\begin{equation}\label{E:1.5'u}
 \vert \nabla ( \Delta^{-1} {\rm div} \, \vec b) 
 \vert^2 + 
 \vert \nabla (\Delta^{-1} \, q)  \vert^2 \in \mathfrak{M}^{1, \, 2}_+(\R^n). 
\end{equation} 
\end{main1}

\begin{remark} In the case $n=2$,  we will show that  {\rm(\ref{E:1.1a})} holds if 
and only if $\tfrac 12 \, (A+ A^t) \in L^\infty(\R^2)^{2\times 2}$,  
$\Delta^{-1} {\rm curl} \,\vec b - \tfrac 1 2 \, (A- A^t) 
\in {\rm BMO} (\R^2)^{2\times 2}$, and $q=0$.
\end{remark}

\begin{remark} Expressions like 
$\nabla ( \Delta^{-1} {\rm div} \, \vec b)$, ${\rm Div} 
( \Delta^{-1} {\rm curl} \, \vec b)$, and $\nabla (\Delta^{-1} \, q)$   
used above 
which involve nonlocal operators   
are defined in the sense of distributions. This is possible, as we 
demonstrate below, since 
 $\Delta^{-1} {\rm div} \, \vec b$, 
$\Delta^{-1} {\rm curl} \, \vec b$, and  $\Delta^{-1} q$
can be understood in terms of the convergence in the weak-$*$ topology 
of ${\rm BMO}(\R^n)$    of, respectively, $\Delta^{-1} \, {\rm div} \, (\psi_N \, \vec b)$, 
$\Delta^{-1} \, {\rm curl} \, (\psi_N \, \vec b)$, and 
$\Delta^{-1} \,  (\psi_N \, q)$ as $N \to +\infty$. 
Here  $\psi_N$ is a smooth 
cut-off function supported on $\{x: \,  |x|<N\}$, and the limits 
above do not depend on the choice of $\psi_N$.   
\end{remark}

It follows from Theorem~I that $\mathcal L$  is form bounded on $L^{1, \, 2}(\R^n)\times L^{1, \, 2}(\R^n)$ if and 
only if the symmetric part of $A$ is essentially bounded, i.e., 
$\tfrac 1 2 (A + A^t ) \in L^\infty(\R^n)^{n \times n}$, and $\vec b_1 \cdot \nabla +q$ 
is form bounded, where 
  \begin{equation}\label{E:1.1d}
\vec {b_1} = \vec b  - \tfrac 1 2 \, {\rm Div} (A-A^t) . 
\end{equation}

In particular,  the principal part 
$\mathcal P u={\rm div} (A \, \nabla u)$ 
is form bounded 
if and only if  
\begin{align}\label{E:1.1e}
& \tfrac 1 2 (A + A^t)  \in L^\infty(\R^n)^{n \times n}, \\
& \Delta^{-1} [{\rm curl} \, {\rm Div} \,\tfrac 1 2 (A - A^t )]  
\in {\rm BMO}(\R^n)^{n \times n}.
\label{E:1.1ee}
\end{align} 
 A simpler condition with  $\tfrac 1 2 (A-A^t) \in {\rm BMO}(\R^n)^{n \times n}$ 
in place of 
(\ref{E:1.1ee}) is sufficient, but generally not necessary, unless $n\le 2$.

Thus,  the form boundedness problem for the general second 
order differential operator in the divergence form (\ref{E:1.00}) 
is reduced  to the special case 
\begin{equation} \label{E:1.l}
\mathcal L = \vec b\cdot \nabla  + q, \qquad \vec b \in D'(\R^n)^n, \quad q \in D'(\R^n). 
\end{equation}

  As a  corollary of Theorem~I, 
 we deduce that, if $\vec b \cdot \nabla + q$  is form bounded, i.e., for all $u, \, v \in C^\infty_0(\R^n)$, 
\begin{equation}\label{E:1.1100}
\left \vert \int_{\R^n}   ( \vec b\cdot\nabla u \, \, 
\overline{v} + q  u \,  \overline{v} ) \, dx 
\right \vert  \leq C \, ||u||_{L^{1, \, 2} (\R^n)}  
||v||_{L^{1, \, 2} (\R^n)},     
\end{equation}
then the Hodge decomposition 
\begin{equation}\label{E:1.4uu}
\vec b = \nabla ( \Delta^{-1} {\rm div} \, \vec b) + {\rm Div} \, (\Delta^{-1} {\rm curl} \, \vec b)  
\end{equation} 
holds, where $\Delta^{-1} {\rm curl} \, \vec b \in {\rm BMO}(\R^n)^{n \times n}$, 
and 
\begin{equation}\label{E:1.7u} 
\int_{|x-y|<r}   \,  [ \,    \vert  
\nabla ( \Delta^{-1} {\rm div} \, \vec b)  
\vert^2 +  \vert \nabla (\Delta^{-1} \, q)  \vert^2 \, ] \,   
dy \le {\rm const} \, \,  r^{n-2},  
\qquad 
\end{equation}
for  all 
 $r>0, \, x \in \R^n$, in the case $n\ge 3$; in two dimensions, it follows that ${\rm div}\, \vec b = q =0$.  

We observe that   condition (\ref{E:1.7u})  is generally stronger than 
 $\Delta^{-1} {\rm div} \, \vec b \in 
{\rm BMO}$ and $\Delta^{-1} \, q \in {\rm BMO}$, while  the divergence-free 
part of $\vec b$ is characterized by   
$\Delta^{-1} {\rm curl} \, \vec b \in {\rm BMO}$, for all $n \ge 2$.

A close sufficient condition of the  Fefferman--Phong type  
can be stated in the following form: 
\begin{equation}\label{E:1.8u} 
\int_{|x-y|<r}   \, [ \,  \vert  
\nabla ( \Delta^{-1} {\rm div} \, \vec b)  
\vert^2 + \vert \nabla (\Delta^{-1} \, q)   \vert^2
\,  ]^{1 + \epsilon} 
\,   dy \le {\rm const} \, r^{n-2(1 + \epsilon) }, 
\end{equation} 
for some $\epsilon > 0$ and all $r>0, \, x \in \R^n$. This is a consequence of  Theorem~I 
coupled  with  (\ref{E:1.5b}),  
where  $\vert ( \Delta^{-1} {\rm div} \, \vec b)  
\vert^2 + \vert \nabla (\Delta^{-1} \, q)   \vert^2$ is used in place of $q$. Sharper conditions of the Chang--Wilson--Wolff  type are readily deduced from Theorem I  
 by combining it with the results of  \cite{ChWW}.

 It is worth mentioning 
that the class 
of potentials obeying  (\ref{E:1.8u})  is substantially broader than its  
subclass 
\begin{equation}\label{E:1.5aa} 
\int_{|x-y|<r}   ( | \vec b|^2 + |q| )^{1 + \epsilon} 
\,   dy \le {\rm const} \, r^{n-2(1 + \epsilon) }.  
\end{equation} 
The sufficiency of the preceding condition for (\ref{E:1.1100}) 
is deduced  by a direct application of the original Fefferman--Phong condition  and Schwarz's inequality.

More generally,  (\ref{E:1.1100})   
clearly  follows from a cruder estimate, 
\begin{equation}\label{E:1.2''}
\int_{\R^n}  |u|^2  \,  ( |\vec b|^2 + |q|  ) \, 
dx \le {\rm const} \, 
|| u||^2_{L^{1, \, 2}(\R^n)}, 
\quad   u \in C^\infty_0(\R^n),
\end{equation}
which is equivalent to $|\vec b|^2 + |q| \in \mathfrak{M}^{1, \, 2}_+(\R^n)$. 

However, by replacing (\ref{E:1.1100})  with (\ref{E:1.2''}),  
one strongly reduces the class of admissible vector fields $\vec b$ 
and potentials $q$. Various examples of this phenomenon in the case $\vec b = 0$ 
are given in \cite{MV1}. 
An instructive example for $\vec b \cdot \nabla$ in the case $q=0$ is provided by the vector field 
$$
\vec b (x) =  \left ( {x_2} (x_1^2 + x_2^2)^{-1}, 
- {x_1} (x_1^2 + x_2^2)^{-1}, 0, \ldots 0 \right ), 
\quad x \in \R^n, 
$$
where $n \ge 2$. An elementary argument involving 
polar coordinates 
and a Fourier series expansion shows that this vector field obeys 
(\ref{E:1.1100}). On the other hand, (\ref{E:1.2''}) fails since 
$\vec b \not\in L^2_{{\rm loc}}(\R^n)$. 

We note in passing that, for $q=0$,  (\ref{E:1.2''}) is equivalent 
to the boundedness of 
the nonlinear quadratic form $\langle \,  | \vec b \cdot 
\nabla u|, \, u \rangle$ on $L^{1, \, 2} (\R^n)\times L^{1, \, 2} (\R^n)$ 
(see Sec.~\ref{Section 4}). As it turns out, dealing with the linear version 
$\langle  \vec b \cdot \nabla u, \, u \rangle  $ is more difficult.

The main obstacle in the proof of Theorem I  is the interaction between 
the quadratic forms associated with  $q - \tfrac 1 2 \, {\rm div} \, \vec b $ and 
the divergence free part of $\vec b$ (see Sections~\ref{Section 2a} and \ref{Section 3}). To overcome this difficulty, one needs 
to distinguish  the class of vector fields $\vec b$ such that the commutator 
inequality 
\begin{equation}\label{E:1.1aa}
\left \vert \int_{\R^n} \vec b \cdot (u \, \nabla \bar v - \bar v \, 
\nabla u) \,  \, dx 
\right \vert  \leq \, {\rm const} \, ||u||_{L^{1, \, 2} (\R^n)} \, 
||v||_{L^{1, \, 2} (\R^n)}    
\end{equation}
holds for all $u, v  \in C^\infty_0(\R^n)$. In the important special 
case of irrotational fields where 
$\vec b = \nabla f$, the preceding inequality is equivalent to
the boundedness of the commutator   $[ f, \, \Delta]$ 
acting from
$L^{1, \, 2}(\R^n)$ to  $L^{-1, \, 2}(\R^n)$.

A complete characterization of those $\vec b$ which obey (\ref{E:1.1aa}) 
is obtained below (Sec.~\ref{Section 3}, Lemma~\ref{Lemma 3.6}) using the idea of 
the gauge transformation  (\cite{LL}, Sec. 7.19; \cite{RS}, Sec. X.4): 
$$\nabla  \to   e^{- i \lambda} \, \nabla \, e^{+ i \lambda},
$$
where the gauge $\lambda$ is 
a real-valued  function which lies in $L^{1, \, 2}_{{\rm loc}}(\R^n)$. 

The problem of choosing an appropriate gauge is known to be highly nontrivial. 
In the present paper, $\lambda$ is picked in a very specific form:
$$\lambda = \tau \, \log \, (P \mu), \qquad 1< 2 \tau < 
\tfrac n {n-2}, \quad n \ge 3,$$
where $\tau$ is a constant, and $P \mu = (-\Delta)^{-1} \mu$ 
is the Newtonian potential of the equilibrium measure $\mu$
associated with an arbitrary compact set $e$ of positive  capacity, 
 ${\rm cap} \, (e)>0$ (see the definitions in Sec.~\ref{Section 2}). 

We will verify
 that, with this choice of $\lambda$,  the energy space $L^{1, \, 2} (\R^n)$ is  gauge 
invariant, and the irrotational part 
  $\vec c= \nabla ( \Delta^{-1} {\rm div} \, \vec b)$ of $\vec b$   
obeys 
$$\int_e |\vec c |^2 \, dx \le {\rm const} \, {\rm cap} \, (e),$$
where the constant does not depend on $e$. This is known 
to be equivalent 
to $|\vec c|^2  \in \mathfrak{M}^{1, \, 2}_+(\R^n)$ (see Theorem~\ref{Trace Theorem}). 
In addition, a careful analysis 
shows that
 $F=\Delta^{-1} {\rm curl} \, \vec b$  
belongs to ${\rm BMO}$, and $\vec b = \vec c + {\rm Div} \, F$.
  These conditions 
combined turn out to be necessary and sufficient for (\ref{E:1.1aa}).

At the end of  Sec.~\ref{Section 3}, we give applications to the 
magnetic Schr\"odinger operator $\mathcal{M}$  
defined by (\ref{E:1.1m}). We show 
that $\mathcal{M}$ is 
form bounded if and only if  both 
$q + | \vec a|^2$ and $\vec a \cdot \nabla$ are form bounded. Thus, the  
form boundedness criterion of $\mathcal{M}$ is deduced from Theorem I (see Theorem~\ref{Theorem 3.4}).

In Sec.~\ref{Section 5}, we extend our results to 
the Sobolev space 
$W^{1, \, 2}(\R^n)$. 
In particular, we give necessary and sufficient 
conditions  (Theorem~\ref{Theorem 5.1})  
for the boundedness of the  general second order operator 
$$\mathcal{L} \, : \,  W^{1, \, 2}(\R^n) \to  W^{-1, \, 2}(\R^n).$$
 This solves the {\it relative form boundedness} problem  for $\mathcal{L}$,  and consequently for 
the magnetic Schr\"odinger operator $\mathcal{M}$,  
 with respect to the Laplacian  on $L^2(\R^n)$ (see \cite{RS}, Sec. X.2). The proofs are based on a localized 
version of the approach developed in Sec.~\ref{Section 3}, and 
in particular involve an 
 inhomogeneous version of the ${\rm div}$-${\rm curl}$ lemma 
(see Lemma~\ref{Lemma 5.2} below).

We remark that other fundamental properties of quadratic  
  forms associated with  differential operators  can be characterized using our methods.  
 For the Schr\"odinger operator 
$\mathcal{L} = \Delta + q$, criteria of relative  compactness 
were obtained in \cite{MV1}, while the infinitesimal form boundedness expressed by the inequality \begin{equation}\label{E:1.9} 
| \langle \mathcal{L} \, u, \, u \rangle | \le 
\epsilon \, ||\nabla u||^2_{L^2(\R^n)} + 
C(\epsilon) \, || u||^2_{L^2(\R^n)}, \quad   u \in C^\infty_0 
(\R^n), 
\end{equation} 
for every $\epsilon > 0$, along with  Trudinger's condition where 
$C(\epsilon) = C \, \epsilon^{-\beta}$, $\beta>0$,  was characterized in \cite{MV4}. 
 Necessary and sufficient conditions for such properties in the case of the general second order differential operator  
are discussed in Sec.~\ref{Section 6}.

\section{Preliminaries}\label{Section 2}

By $L^{1, \, 2}(\R^n)$ we denote the energy space (homogeneous Sobolev space)
 defined in the Introduction as the completion of the complex-valued   
$C^\infty_0$ functions in the Dirichlet norm. For $n \ge 3$, an
 equivalent norm on $L^{1,2}(\R^n)$ is given 
by 
$$ 
||u||_{L^{1,2}(\R^n)} = ||\,  |x|^{-1} u  ||_{L^2(\R^n)} + 
||\nabla u||_{L^2(\R^n)}, 
\qquad u \in L^{1,2}(\R^n).    
$$ 

By 
$W^{1,2}(\R^n)$ we denote the space of 
weakly differentiable (complex-valued) functions on $\R^n$ $(n\ge 1)$ such that 
$$
||u||_{W^{1,2}(\R^n)} = ||u||_{L^2(\R^n)} + ||\nabla u||_{L^2(\R^n)}< +\infty.  
$$
The dual spaces are respectively $L^{-1,2} (\R^n) = L^{1,2} (\R^n)^*$ and   
 $W^{-1,2} (\R^n) = W^{1,2} (\R^n)^*$.

For $0< r < \infty$, denote  by $L^r_{\rm unif}(\R^n)$ 
all $f\in L^r_{\rm loc}(\R^n)$ such that 
$$
||f||_{L^r_{\rm unif}} = \sup_{x_0 \in \R^n} \, 
||\chi_{B_1 (x_0)} \,  f ||_{L^r(\R^n)} < \infty.
$$

 We set 
$$
m_B (f) = \frac 1 {|B|} \int_B f(x) \, dx
$$
for a ball $B\subset \R^n$, and denote by  
 ${\rm BMO}(\R^n)$ the class of 
$f \in L^r_{\rm loc} (\R^n)$ for which 
$$
\sup_{x_0\in \R^n, \, \delta>0} \, \,  \frac {1}{|B_\delta(x_0)|} 
\int_{B_\delta (x_0)} |f(x)-m_{B_\delta(x_0)}(f)|^r \, dx < + \infty,
$$
for any (or, equivalently, all) $1\le r < +\infty$. 
An inhomogeneous version of ${\rm BMO}(\R^n)$ (the so-called 
local ${\rm BMO}$; see  \cite{St}, p. 264), which we denote by 
${\rm bmo}  (\R^n)$, is  defined  as 
the set of  $f \in L^r_{\rm unif}(\R^n)$ such that the preceding condition holds 
for all $0<\delta\le 1$, and additionally 
$$\sup_{x_0\in \R^n, \, \delta \ge 1} \, \, \frac {1}{|B_\delta(x_0)|} 
\int_{B_\delta (x_0)} |f(x)|^r \, dx < + \infty.
$$

We will also need the space  ${\rm BMO}^\#(\R^n)$  defined as the set 
of $f \in L^r_{{\rm loc}}(\R^n)$ 
such that 
$$
\sup_{x_0 \in \R^n, \, 0<\delta\le 1} \, \frac 1 {|B_\delta(x_0)|} 
\int_{B_\delta(x_0)} |f(x) - m_{B_\delta(x_0)}(f)|^r \, dx < +\infty,
$$
for any (or equivalently all) $1 \le r < +\infty$. 
Notice that ${\rm bmo(\R^n)} \subset {\rm BMO}(\R^n)\subset 
{\rm BMO}^\#(\R^n)$.

The corresponding vector- and matrix-valued function spaces are introduced in 
a similar way. In particular,  ${\rm BMO}(\R^n)^n$ stands for 
 the class of vector fields 
$\vec f = \{f_j\}_{j=1}^n : \, \R^n \to \C^n$, such that 
$f_j \in {\rm BMO}(\R^n)$, $j =1, 2, \dots, n$. The 
 matrix-valued analogue is denoted by ${\rm BMO}(\R^n)^{n \times n}$, etc.

 For a matrix field $F = (f_{i j})_{i, j =1}^n \in D'(\R^n)^{n \times n}$, 
the matrix divergence 
operator  
${\rm Div}$ is defined by ${\rm Div} \, F = \left ( \sum_{j=1}^n \, \partial_j \, 
f_{i j}  \right)_{i=1}^n \in D'(\R^n)^n$. The Jacobian, $\D$, is 
the formal adjoint of $-{\rm Div}$ (see, e.g., \cite{IM}) :
$$\langle  {\rm Div} \, F,  \, \vec v \rangle = - {\rm trace} \, 
\langle F^t, \, \D \, \vec v 
\rangle, \qquad   \, \vec v \in C^\infty_0(\R^n)^n,$$
where $F^t = (f_{ji})_{i, j =1}^n$ is the transposed matrix field. If $F$ is 
skew-symmetric, i.e., $f_{ij} = -f_{ji}$, 
then obviously ${\rm div} \, ({\rm Div} \, F) = 0$.

The capacity of a compact set $e\subset \R^n$ is defined by 
(\cite{LL}, Sec. 11.15; \cite{M}, Sec. 2.2): 
\begin{equation}\label{E:cap}
{\rm cap} \, (e) = \inf \, 
\left\{ \, ||u||^2_{L^{1,2} (\R^n)}  : \quad 
u \in C^\infty_0(\R^n), \quad u(x) \ge 1 \, \,  {\rm on} 
\, \,  e \right\}.
\end{equation}
For a cube or ball $Q$ in $\R^n$, 
\begin{equation}\label{E:cube}
{\rm cap} \, (Q) \simeq |Q|^{1-\frac 2 n} \quad {\rm if} \, \,  n \ge 3; 
\quad {\rm cap} \, (Q) =0 \quad  {\rm if} \, \, n=2.
\end{equation}

We will also need the capacity ${\rm Cap} \, (\cdot)$ associated with the 
Sobolev space $W^{1, \, 2}(\R^n)$ defined  by
\begin{equation}\label{E:Cap}
{\rm Cap} \, (e) = \inf \, 
\left\{ \, ||u||^2_{W^{1,2} (\R^n)}  : \quad 
u \in C^\infty_0(\R^n), \quad u(x) \ge 1 \, \,  {\rm on} 
\, \,  e \right\},
\end{equation}
for compact sets $e\subset \R^n$. 
Note that ${\rm Cap} \, (e) \simeq {\rm cap} \, (e)$ if ${\rm diam} \, (e) \le 1$, and $n \ge 3$.  
For a cube or ball $Q$ in $\R^n$,  
\begin{equation}\label{E:Cube}
{\rm Cap} \, (Q) \simeq |Q|^{1-\frac 2 n} \quad {\rm if} \, \,  n \ge 3; 
\quad {\rm Cap} \, (Q) \simeq \left (\log \tfrac 2 {|Q|}\right)^{-1} \quad  {\rm if} \, \, n=2,
\end{equation}
provided $|Q| \le 1$. For these and other properties of capacities, as well as related notions of 
potential theory we refer to \cite{AH},  \cite{M}.

We conclude this section with several equivalent characterizations of the class of admissible 
measures $\mu \in \mathfrak{M}_+^{1, \,2}(\R^n)$ which obey the trace inequality 
\begin{equation}\label{E:tr0}
 \int_{\R^n} |u|^2 \,  d \mu \le c^2 \, ||u||^2_{L^{1, \, 2}(\R^n)}, \qquad u \in C^\infty_0(\R^n),
\end{equation}
where $c$ is a positive constant which does not depend on $u$. 

By 
$(-\Delta)^{-\frac 1 2} \mu = c(n) \, \int_{\R^n} |x-t|^{1-n} \, d \mu(t)$ we denote the Riesz potential of order $1$ 
of the measure $\mu$; here $c(n)$ is a normalization constant which depends only on $n$.

\begin{theorem}\label{Trace Theorem} Let $\mu$ be a locally finite nonnegative measure 
on $\R^n$. Then  $\mu \in  \mathfrak{M}_+^{1, \,2}(\R^n)$ if and only if any one of the following 
statements hold.

{\rm (i)}  The Riesz potential $(-\Delta)^{-\frac 1 2} \mu \in L^2_{{\rm loc}} (\R^n)$, and 
$[(-\Delta)^{-\frac 1 2} \mu]^2 \in \mathfrak{M}_+^{1, \,2}(\R^n)$, i.e., 
 \begin{equation}\label{E:tr1}
 \int_{\R^n} |u|^2 \, [(-\Delta)^{-\frac 1 2} \mu]^2 \,  d x \le c_1^2 \, ||u||^2_{L^{1, \, 2}(\R^n)}, 
 \qquad u \in C^\infty_0(\R^n),
\end{equation}
where $c_1>0$ does not depend on $u$.

{\rm (ii)} For every compact set $e \subset \R^n$, 
 \begin{equation}\label{E:tr2}
 \mu (e) \le c_2 \, \text{\rm{cap}} \, (e),
\end{equation}
where $c_2$ does not depend on $e$.

{\rm (iii)} For every  ball $B$ in $\R^n$,
 \begin{equation}\label{E:tr3}
 \int_B [(-\Delta)^{-\frac 1 2} \mu_B]^2 \, dx \le c_3 \, \mu (B),
\end{equation}
where $d \mu_B = \chi_B \, d \mu$, 
and $c_3$ does not depend on $B$.

{\rm (iv)} The pointwise inequality
 \begin{equation}\label{E:tr4}
(-\Delta)^{-\frac 1 2} [ (-\Delta)^{-\frac 1 2} \mu   ]^2 (x) \le c_4 \, (-\Delta)^{-\frac 1 2} \mu(x) < \infty 
\end{equation}
holds a.e., where $c_4$ does not depend on $x \in \R^n.$ 

{\rm (v)} For every dyadic cube $P$ in  $\R^n$, 
\begin{equation}\label{E:tr5}
\sum_{Q \subseteq P}  \left [ \frac {\mu(Q)}{
|Q|^{1 - 1/n}} \right ]^2 |Q| \le c_5 \, \mu (P),
\end{equation}
where the sum is taken over all dyadic cubes $Q$ contained in $P$, and 
$c_5$ does not depend on $P$.

Moreover, the least constants $c_i$, $i=1, \ldots, 5$,  are equivalent to the least constant $c$ in {\rm (\ref{E:tr0})}. 
\end{theorem}

Theorem~\ref{Trace Theorem} follows from the results of \cite{KS}, \cite{M}, \cite{MV1}, and 
\cite{V}.

 \begin{remark} An analogous characterization holds for admissible measures on the space $W^{1, \, 2}(\R^n)$ 
in place of $L^{1, \, 2}(\R^n)$. One only needs to replace $(-\Delta)^{-\frac 1 2} \mu$ in 
statements {\rm (i)},  {\rm (iii)},  and {\rm (iv)}  
by 
 $(1-\Delta)^{-\frac 1 2} \mu$, the capacity ${\rm cap} \, (\cdot)$ in   {\rm (ii)}
   by ${\rm Cap} \, (\cdot)$, and restrict oneself to cubes $P$ such that $|P| \le 1$ in 
{\rm (v)}. 
\end{remark}

\section{Reduction to inequalities for  lower order terms}\label{Section 2a}

In this section,  the form boundedness problem for the general second order 
differential operator $\mathcal L$  defined by (\ref{E:1.0}) is reduced to the special case of   
 lower order terms, $\vec b \cdot \nabla + q$. 
The latter, in its turn, is shown to be equivalent to the form boundedness of $q - \tfrac 1 2  \, {\rm div} \, \vec b$,
 and  the commutator inequality (\ref{E:1.1aa}). 

Since the coefficients $A=(a_{ij})$,  $\vec b = (b_i)$, and $q$  are arbitrary real- or complex-valued distributions, 
we may assume without loss of generality 
that $\mathcal L$ is in the divergence form, 
\begin{equation}\label{E:2.00}
\mathcal L \, u= {\rm div} \, (A \, \nabla u) + \vec b \cdot \nabla u + q \, u,  \qquad u \in C^\infty_0(\R^n),  
\end{equation}
with the same  principal part as  (\ref{E:1.0}).  We denote by 
$A^{s} = \frac 12 (A +A^t)$ and $A^{c} = \frac 12 (A -A^t)$ respectively 
the symmetric and skew-symmetric parts of $A$, and by ${\rm Div}$ the row divergence operator 
acting from $D'(\R^n)^{n \times n}$ to  $D'(\R^n)^{n}$.

\begin{proposition} \label{Proposition 2.100}  Suppose $\mathcal L$ is defined by  {\rm (\ref{E:2.00})}, 
where $A= (a_{ij})_{i, \, j=1}^n \in D'(\R^n)^{n \times n}$, 
$\vec b = (b_j)_{j=1}^n\in D'(\R^n)^{n}$, and $q \in D'(\R^n)$, $n \ge 1$. Let $\vec b_1 = 
\vec b - {\rm Div} \, A^{c}$. Then the following statements are equivalent. 

{\rm (i)} The sesquilinear form associated with $\mathcal L$ is bounded, i.e., 
\begin{equation}\label{E:2.01}
| \langle \mathcal{L} \, u, \, v \rangle | \leq \, C \, ||u||_{L^{1, \, 2} (\R^n)} \, 
||v||_{L^{1, \, 2} (\R^n)}    
\end{equation}
where the constant $C$ does not depend on $u, \, v  \in C^\infty_0(\R^n)$.

{\rm (ii)} The   sesquilinear form associated with  $\mathcal L_1 = \vec b_1 \cdot \nabla + q$ is bounded, 
and $A^{\rm s} \in L^\infty(\R^n)^{n \times n}$.  
\end{proposition}

\begin{proof}  Since $A^{c}$ is skew-symmetric, 
${\rm div} \, (A^{c} \, \nabla u) = - {\rm Div} \, A^{c} \cdot \nabla u,$ 
and consequently 
$$\langle {\rm div} \, ( A \, \nabla u), \, v \rangle = - \langle A^{s} \, \nabla u, \, \nabla v \rangle
- \langle  {\rm Div} \, A^{c} \cdot \nabla u, \,  v \rangle,
$$
for $u, \, v \in  C^\infty_0(\R^n)$.
Hence, 
$$\langle \mathcal{L} \, u, \, v \rangle = - \langle A^{\rm s} \, \nabla u, \, \nabla v \rangle + 
\langle \vec b_1 \cdot \nabla \, u, \, v \rangle + \langle q \, u, \, v \rangle,
$$
where $\vec b_1 = \vec b - {\rm Div} \, A^{\rm c} $.

Suppose that the sesquilinear form of $\mathcal{L}$ is  bounded. 
Then, replacing $u$ and $v$ in   (\ref{E:2.01})
respectively  by $\tilde u = e^{i t \, \xi \cdot x} \, u$ and 
 $\tilde v = e^{i t \, \xi \cdot x} \, v$, 
where $x, \, \xi \in \R^n$, and $t>0$ , 
we obtain 
\begin{align}
& \vert - \langle A^{\rm s}  \, \nabla \tilde u, \, \nabla \tilde v \rangle  + \langle  \vec b_1 \cdot \nabla \tilde  u, \, \tilde v\rangle 
+  \langle q \, \tilde u, \, \tilde v \rangle 
\vert   \le  C \,  ||\tilde u||_{L^{1, \, 2}(\R^n)} \, ||\tilde v||_{L^{1, \, 2}(\R^n)}
\notag \\  & \le C \, ( ||u||_{L^{1, \, 2}(\R^n)} +    t \, |\xi| \, ||u||_{L^{2}(\R^n)})  \, 
 ( ||v||_{L^{1, \, 2}(\R^n)} +    t \, |\xi| \, ||v||_{L^{2}(\R^n)}).
\notag
\end{align}
Dividing both sides of the preceding inequality by $t^2$, and letting $t \to +\infty$, we see that the last two terms 
on the left-hand side tend to $0$, which yields 
$$\left \vert \int_{\R^n} ( A^{\rm s} \, \xi\cdot \xi ) \, u \, \bar v \, dx  \right \vert \le  C \,  
|\xi|^2 \, ||u||_{L^{2}(\R^n)} \, ||v||_{L^{2}(\R^n)}.
$$
From this we deduce 
$$| A^{s}(x) \, \xi \cdot \xi |  \le C \, |\xi|^2, \qquad x, \, \xi \in \R^n.$$
Clearly, both  the real and imaginary parts of $A^{s}$ obey 
the preceding inequality,  and since $A^{s}$ is symmetric, their operator norms are bounded by $C$. 
Hence, necessarily $A^{s} \in L^\infty(\R^n)^{n \times n}$. The latter is also 
sufficient for the form boundedness of ${\rm div} \, (A^{s} \, \nabla)$. Thus, $\mathcal{L} $ is form bounded if and only if 
 $A^{\rm s} \in L^\infty(\R^n)^{n \times n}$,  and $\vec b_1 \cdot \nabla +q$ is form bounded. \end{proof}

\begin{proposition} \label{Proposition 3.1}
 Let  $q \in D'(\R^n)$ and $\vec b \in D'(\R^n)^n$, $n \ge 1$, and let $\mathcal{L} = \vec b \cdot \nabla + q$. 
Then the following statements are equivalent. 

{\rm (i)} The bilinear form associated with 
$\mathcal{L}$ obeys {\rm (\ref{E:2.01})}.

{\rm (ii)} The following two conditions hold:

{\rm (a)} For all $u, \, v \in C^\infty_0(\R^n)$, 
\begin{equation}\label{E:3.2}
\left \vert 
 \langle (q - \tfrac 1 2 \, {\rm div} \, \vec b) 
\, u, \, v \rangle \right \vert \le C \, 
||u||_{L^{1, \, 2}(\R^n)} \, ||v||_{L^{1, \, 2} (\R^n)}. 
\end{equation} 

{\rm (b)} For all $u, \, v \in C^\infty_0(\R^n)$, 
 \begin{equation}\label{E:3.3}
\left \vert 
 \langle  \vec b,  \, 
 \bar u \, \nabla v  - v \nabla \bar u \rangle \right \vert 
\le C \, 
||u||_{L^{1, \, 2} (\R^n)} \, ||v||_{L^{1, \, 2} (\R^n)}.
\end{equation} 

\end{proposition}

\begin{proof} Integration by parts gives 
$$ \langle \vec b \cdot \nabla u  + q \, u, \, v \rangle  = 
\langle (q - \tfrac 1 2 \, {\rm div} \, \vec b),  \, \bar u \, v \rangle  
- \tfrac 1 2 \langle \vec b, \, \bar u \, \nabla v - v \, \nabla \bar u 
\rangle.
$$
Interchanging the roles of $\bar u$ and $v$, it is easy to see
 that the bilinear form associated with $\vec b \cdot \nabla   + q$ 
is bounded if and only 
if both forms on the right-hand side of the preceding 
equation are bounded, i.e., both 
(\ref{E:3.2}) and (\ref{E:3.3}) hold.
\end{proof}

\begin{remark} Inequality {\rm (\ref{E:3.2})} holds if and only if the 
 inequality 
\begin{equation}\label{E:3.3'}
\int_{\R^n}   \vert \nabla \Delta^{-1} \, 
( q - \tfrac 1 2 \, {\rm div} \, \vec b ) \vert^2 \, |u|^2 \, dx 
 \le C \, 
||\nabla u||^2_{L^2(\R^n)} 
\end{equation} 
is valid, where $C$ does not depend on
 $u \in C^\infty_0(\R^n)$ (see \cite{MV1}, Theorem I). 
\end{remark}

\begin{corollary} \label{Corollary 2.100} Let 
$\mathcal P u= {\rm div} \, (A \nabla u)$, 
where $A= (a_{ij})_{i, \, j=1}^n \in D'(\R^n)^{n \times n}$. Then 
\begin{equation}\label{E:2.010}
| \langle \mathcal{P} \, u, \, v \rangle | \leq \, C \, ||u||_{L^{1, \, 2} (\R^n)} \, 
||v||_{L^{1, \, 2} (\R^n)}    
\end{equation}
for all $u, \, v  \in C^\infty_0(\R^n)$, if and only if $A^s \in L^\infty(\R^n)^{n \times n}$, 
and $A^c$ obeys the inequality
\begin{equation}\label{E:2.011}
| \langle {\rm Div} \, A^c,  \, u \, \nabla v - v \, \nabla u\rangle | \leq \, C \, ||u||_{L^{1, \, 2} (\R^n)} \, 
||v||_{L^{1, \, 2} (\R^n)}    
\end{equation}
for all $u, \, v  \in C^\infty_0(\R^n)$. The preceding inequality holds if 
$A^c \in {\rm BMO} (\R^n)^{n \times n}$. 
\end{corollary} 

Corollary~\ref{Corollary 2.100} follows from 
 Propositions~\ref{Proposition 2.100} 
and \ref{Proposition 3.1}. The last statement is a consequence of 
the div-curl lemma \cite{CLMS} (see \cite{T}, Sec. 3.8). 
A more precise necessary and sufficient condition for (\ref{E:2.011}) is obtained below. 

\section{Operators with lower order terms and magnetic 
Schr\"odinger operators}\label{Section 3}

In this section, which contains our main results, 
 we  consider the form boundedness problem on 
$L^{1, \, 2} (\R^n)$, $n \ge 2$,  for 
the  operator 
\begin{equation}\label{E:3.1}
\mathcal{L} =  \vec b \, \cdot \nabla + q
\end{equation}
with distributional lower order terms $\vec b \in D'(\R^n)^n$ 
and $q \in D'(\R^n)$. Here $\mathcal{L}$ is initially defined as 
acting from $D(\R^n)$ to $D'(\R^n)$. We deduce necessary and 
sufficient conditions for the form boundedness of 
$\mathcal{L}$, and as a consequence, of  
the Schr\"odinger operator $(i \, \nabla +\vec a )^2 + q$
with magnetic vector potential $\vec a \in L^2_{{\rm loc}}(\R^n)$. 

We will need a series of  lemmas and propositions.

\begin{proposition} \label{Proposition 3.2}
 Let   $\vec b \in D'(\R^n)^n$, $n \ge 2$. Suppose that 
 {\rm (\ref{E:3.3})} holds. 
Then, for every cube $Q$ 
in $\R^n$, the following estimates are valid:
\begin{align}\label{E:3.4p}
|| {\rm div} \, \vec b||_{L^{-1, \, 2}(Q)} & \le C \,
|Q|^{\frac 1 2 - \frac 1 n} \quad {\rm if} \, \, n \ge 3;
\qquad {\rm div} \, \vec b = 0 \quad {\rm if} \, \, n =2,
\\ ||\vec b||_{L^{-1, \, 2}(Q)} & \le C \, |Q|^{\frac 1 2}
\quad {\rm if} \, \, n \ge 2,
\label{E:3.4q}
\end{align}
where $C$ does not depend on $Q$. \end{proposition}

\begin{proof} Let $v \in C^\infty_0(Q)$, and suppose 
$u=1$ on $Q$, $u \in C^\infty_0(\R^n)$ in (\ref{E:3.3}). Then   
$$\left \vert \langle  \vec b,  \, \bar u \, \nabla v - v \, \nabla \bar u 
\rangle \right \vert =  \left \vert \langle  \vec b,  \, \nabla v  
\rangle \right \vert = \left \vert \langle  {\rm div} \, \vec b,  \, v  
\rangle \right \vert
\le C \, 
||\nabla u||_{L^2(\R^n)} \, ||\nabla v||_{L^2(Q)}.
$$
Taking the infimum over all such $u$ on the right-hand side, we obtain
$$
\left \vert \langle  {\rm div} \, \vec b,  \, v  
\rangle \right \vert
\le C \, {\rm cap} \, (Q)^{\frac 1 2} \, ||\nabla v||_{L^2(Q)}, 
\qquad   v \in C^\infty_0(Q),
$$
where the capacity ${\rm cap} \, (\cdot)$ is defined by (\ref{E:cap}). 
Taking into account (\ref{E:cube}), 
we deduce from the preceding inequality  that 
${\rm div} \, \vec b =0$ if $n=2$, and 
\begin{equation}\label{E:3.5}
\left \vert \langle  {\rm div} \, \vec b,  \, v 
\rangle \right \vert \le C \, |Q|^{\frac 1 2  - \frac 1 n}
  \, ||\nabla v||_{L^2(Q)}, \qquad   v \in C^\infty_0(Q),  
\end{equation}
if $n=3$, which proves (\ref{E:3.4p}).

Now suppose $v \in C^\infty_0(Q)$, and
 let us set $u = (x_i-a_i) \, \eta \,$  ($i = 1, \ldots, n$), 
where $a=(a_i)$ is the center 
of $Q$, $\eta = 1$ on $Q$ and $\eta \in C^\infty_0(2Q)$. 
Then it is easy to see that 
$||\nabla u||_{L^2(2Q)} \le C \, |Q|^{\frac 1 2}$. Notice that for 
such $u$ and $v$, 
\begin{align}
\langle  \vec b, \, \bar u \, \nabla v - v \, \nabla \bar u 
\rangle  & =    - \langle  {\rm div} \,  \vec b, \, 
 \bar u \,  v \rangle - 2 \, \langle  \vec b, \, v \, 
\nabla \bar u \rangle \notag \\ & = 
- \langle  {\rm div} \,  \vec b, 
 (x_i-a_i) \,  v \rangle - 2 \langle b_i, \, v \rangle.
\notag
\end{align}
Using (\ref{E:3.5}) with $(x_i-a_i) \,  v$ in place of $v$, 
and Poincar\'e's inequality, we obtain
\begin{align}
\left \vert \langle  {\rm div} \,  \vec b, 
 (x_i-a_i) \,  v \rangle \right \vert & \le C \, 
|Q|^{\frac 1 2 - \frac 1 n} \, ||\nabla [(x_i-a_i) v]||_{L^2(Q)} 
\notag \\& \le C \, |Q|^{\frac 1 2 - \frac 1 n} \, ( ||v||_{L^2(Q)} 
+ ||(x_i-a_i) \nabla v||_{L^2(Q)}) \notag \\ 
& \le C \, |Q|^{\frac 1 2} \, ||\nabla v||_{L^2(Q)},
\notag
 \end{align}
for every $v \in C^\infty_0(Q)$. 
Hence,  for every $i=1, \ldots, n$, 
\begin{align}
& 2 \left \vert \langle  b_i, \,  v  \rangle \right \vert  \le   
 \left \vert \langle  \vec b,  \, 
\bar u \, \nabla v - v \, \nabla \bar u \rangle \right \vert 
+ \left \vert \langle  {\rm div} \,  \vec b, 
 (x_i-a_i) \,  v \rangle \right \vert \notag \\ 
& \le C \, 
 ||\nabla u||_{L^2(2 Q)} \, ||\nabla v||_{L^2(Q)} 
 + C \, |Q|^{\frac 1 2} \, ||\nabla v||_{L^2(Q)}
\notag \\ 
& \le C \,  |Q|^{\frac 1 2} \, ||\nabla v||_{L^2(Q)}.
\notag
\end{align}
This yields $||\vec b||_{L^{-1, \, 2}(Q)} \le C \, 
|Q|^{\frac 1 2}$, which completes the proof of 
Proposition~\ref{Proposition 3.2}.
 \end{proof}

For a fixed cube $Q$ in $\R^n$,  we denote by $\{\eta_j\}_{j=0}^\infty$ 
a smooth partition of unity associated with $Q$; i.e., 
$\eta_0\in C^\infty_0 (2 Q)$, 
$\eta_j \in C^\infty_0(2^{j+1} Q \setminus 2^{j-1} Q)$, $j=1, 2, \ldots$, 
so that 
\begin{align}\label{E:res1} 
0 \le \eta_j (x)\le 1, & \qquad  |\nabla \eta_j(x)| \le C \, (2^j \ell(Q))^{-1},
\quad j = 0, 1, \ldots, \\ & \sum_{j=0}^{\infty} \, \eta_j (x) = 1,  \qquad x \in 
\R^n,  
\label{E:res2} 
\end{align} 
where $\ell(Q)$ denotes the side length of $Q$, and $C$ depends only on $n$.

We will need the following proposition.

\begin{proposition}\label{Proposition 2.2} Let $Q$ be a cube in $\R^n$,
$n \ge 2$.  Let $\{\eta_j\}_{j=0}^{+\infty}$ be
the partition of unity associated with $Q$ defined by 
above. Then the
following estimates hold.

{\rm (i)} For any $v \in C^\infty_0(Q)$ and $j=0, 1, \ldots$,
\begin{equation}\label{E:2.20}
|| \nabla(\eta_j  \partial_i  \partial_m  \Delta^{-1}  v)||_{L^2(2^{j+1} Q)}
\le C  2^{-j (1+\frac n 2)}  ||\nabla v||_{L^2(Q)},
\quad i, m=1, \ldots, n,
\end{equation}
where $C$ depends only on $n$.

{\rm (ii)} For any $v \in C^\infty_0(Q)$ and $j=0, 1, \ldots$,
\begin{equation}\label{E:2.20b}
|| \nabla
( \eta_j \, \partial_i \,  \Delta^{-1} \, v)||_{L^2(2^{j+1} Q)} \le C
\,  2^{-j \frac n 2} \, ||v||_{L^2(Q)},
\quad i =1, \ldots, n,
\end{equation}
where $C$ depends only on $n$.

{\rm (iii)} For any $v \in C^\infty_0(Q)$ such that $\int_Q v(x) \, dx =0$,
and $j=2, 3, \ldots$,
\begin{equation}\label{E:2.20c}
|| \nabla ( \eta_j \, \partial_i \,  \Delta^{-1} \, v)||_{L^2(2^{j+1} Q)}
\le C \, 2^{-j (1 +\frac n 2) } \, |Q|^{-\frac 1 2} \, 
|| v||_{L^1 (Q)}, \quad i=1, \ldots, n,
\end{equation} 
where $C$ depends only on $n$.

{\rm (iv)} Let $n \ge 3$. For any $v \in C^\infty_0(Q)$
such that $\int_Q v(x) \, dx =0$, and $j=2, 3, \ldots$,
\begin{equation}\label{E:2.20d}
|| \nabla ( \eta_j \, \Delta^{-1} \, v)||_{L^2(2^{j+1} Q)} \le C
\, 2^{-j \frac n 2} \, |Q|^{\frac 1 n - \frac 1 2} \,
|| v||_{L^1 (Q)}, \quad i=1, \ldots, n,
\end{equation}
where $C$ depends only on $n$.
\end{proposition}

%%%%%%%%%%%%%%%%%%%%%%%%%%%%%%%%%%%%%%%%%%%%%%%%%%%%%%%%%%%%
%{\rm (v)} For any $v \in C^\infty_0(Q)$ and $j=2, 3, \ldots$,
%\begin{equation}\label{E:2.20d1}
%|| \nabla
%( \eta_j \, \partial_i \, \Delta^{-1} \, v)||_{L^2(2^{j+1} Q)} \le C
%\, 2^{-j \frac n 2} \, |Q|^{- \frac 1 2} \,
%|| v||_{L^1 (Q)}, \quad i=1, \ldots, n,
%\end{equation}
%where $C$ depends only on $n$.
%%%%%%%%%%%%%%%%%%%%%%%%%%%%%%%%%%%%%%%%%%%%%%%%%%%%%%%%%%%%

\begin{proof} Let $v \in C^\infty_0(Q)$. Let $a=a_Q$ be 
the center of $Q$, and  $r=\ell(Q)$ its side length.
We denote by $R_i$ the Riesz transforms, and by 
$R_i R_m$, $i, m = 1, \ldots, n$,  the second order Riesz transforms
on $\R^n$ (see \cite{St}. 
For $j=0, \, 1$, (\ref{E:2.20}) follows from the
boundedness of  $R_i R_m$ on $L^2(\R^n)$,
and Poincar\'e's inequality:
\begin{align}
 || \nabla
( \eta_j \, \partial_i \, \partial_m \, \Delta^{-1} \, v)||_{L^2(2^{j+1} Q)} 
& \le || \nabla \eta_j \, ( \partial_i \, \partial_m \, 
\Delta^{-1} \, v) ||_{L^2(2^{j+1} Q)}
\notag \\ & +
|| \eta_j \, \partial_i \, \partial_m \, \Delta^{-1} \,
 \nabla v||_{L^2(2^{j+1} Q)}
\notag \\ & \le C \, \left ( r^{-1} \,
|| R_i \, R_m \, v||_{L^2(\R^n)} + || R_i \, R_m \, \nabla v||_{L^2(\R^n)}
\right) \notag \\ & \le C \, ( r^{-1} \, || v||_{L^2(Q)} + 
||\nabla v||_{L^2(Q)}) \le C \, ||\nabla v||_{L^2(Q)}.\notag
\end{align}
  For $j\ge 2$, one needs   estimates of the kernels of the operators
$\partial_i \, \Delta^{-1}$ and
$\partial_i \, \partial_m \, \Delta^{-1} = -R_i \, R_m$ which 
are given respectively, up to a constant multiple,  by
$$K^i (x-t) =  \frac{(x_i-t_i)}{|x-t|^n}, \quad 
K^{i, \, m}(x-t) =  \frac {\delta_{im} \, |x-t|^2 - n \,
(x_i-t_i) \, (x_m -t_i)} 
{|x-t|^{n+2}}. $$
 Clearly,  
\begin{align}
|K^i (x-t) - K^i (x-a)| & \le C(n) \frac{|t-a|}{|x-t|^n}, \label{E:esta} \\
|K^{i, \, m} (x-t) - K^{i, \, m} (x-a)| & \le C(n) \frac{|t-a|}{|x-t|^{n+1}},
\label{E:estb}
\end{align}
 if $|t-a|<R$, $|x-t|> 2R$. Using  the preceding estimates with $R = 
c(n) \, 2^{j} r$,
we see that,  for $x\in 2^{j+1} Q \setminus 2^{j-1} Q$:
\begin{align}
 & | \partial_i \, \partial_m \, \Delta^{-1} \, v (x)|  =
\left \vert
 \int_{Q} \left (K^i (x-t) - K^i (x-a) \right ) \, \partial_m v (t) \, dt 
\right \vert \notag \\ & \le \, \int_{Q} |K^i (x-t) - K^i (x-a)| \,
|\nabla v(t)| \, dt
\le C \, r^{1-n} \, 2^{-j n} \, ||\nabla  v||_{L^1(Q)}, \notag \\
&  | \nabla \partial_i \, \partial_m \, \Delta^{-1} \, v (x)|  =
\left \vert
 \int_{Q} \left (K^{i, \, m} (x-t) - K^{i, \, m} (x-a) \right ) \, \partial_m \,
\nabla v (t) \, dt\right \vert \notag \\ 
& \le \, \int_{Q} |K^{i, \, m} (x-t) - K^{i, \, m} (x-a)| \, |\nabla v(t)| \, dt
\le C \, r^{-n} \, 2^{-j (n+1)} \, ||\nabla v||_{L^1(Q)}. \notag
\end{align}
Hence,
\begin{align} 
 || \nabla
( \eta_j \, \partial_i \, \partial_m \,
\Delta^{-1} \, v)||_{L^2(2^{j+1} Q)} & \le
|| \nabla  \eta_j \, (\partial_i \,
\partial_m \, \Delta^{-1} \, v)||_{L^2(2^{j+1} Q)}
 \notag \\
  +  ||\eta \,  \partial_i \, \partial_m \, \Delta^{-1} \, \nabla
v||_{L^2(2^{j+1} Q)} & \le
C \, r^{-\frac n 2} \, 2^{-j (1+\frac n 2)} \, ||\nabla v||_{L^1(Q)} 
\notag \\
& \le C \,
2^{-j (1+\frac n 2)} \, ||\nabla v||_{L^2(Q)},\notag
\end{align}
which gives  (\ref{E:2.20}). 

To prove (\ref{E:2.20b}), notice that for $j=0, \, 1$, it follows
\begin{align}
& || \nabla
( \eta_j \, \partial_i \, \, \Delta^{-1} \, v)||_{L^2(2^{j+1} Q)}   \le
|| \nabla  \eta_j \, (\partial_i \,  \Delta^{-1} \, v) ||_{L^2(2^{j+1} Q)}
\notag \\ &  +
|| \eta_j \, \partial_i \,  \, \Delta^{-1} \, \nabla v||_{L^2(2^{j+1} Q)}
\le C \,  ( r^{-1} \,
|| \nabla \Delta^{- 1} \, v||_{L^2(\R^n)} + \sum_{m=1}^n \,
|| R_i \, R_m \, v||_{L^2(\R^n)} ) \notag \\ & \le C \,
|| \nabla \Delta^{- 1} \, v||_{L^q(Q)}) + C \, ||v||_{L^2(Q)},\notag
\end{align}
where $q=\frac {2n}{n-2}$. Estimating the first term
on the right by means of Sobolev's inequality, we conclude
that it is bounded by $C \, ||v||_{L^2(Q)}$.

If $j=2,3, \dots$, then
for $x\in 2^{j+1} Q \setminus 2^{j-1} Q$ we have:
\begin{align}
& | \nabla
( \eta_j \, \partial_i \, \, \Delta^{-1} \, v)(x)| \le
 | \nabla \eta_j (x)| \, | \partial_i \, \Delta^{-1} \, v (x)|
+ |\eta_j(x)| \, |\nabla  \partial_i \, \Delta^{-1} \, v (x)|
 \notag \\ & \le \, C \, (2^j r)^{-1} \,
 \int_{Q} |K^i (x-t)| \, |v (t)| \, dt +
\sum_{m=1}^n \,
 \int_{Q} | K^{i, \, m} (x-t)| \, |v (t)| \, dt
\notag \\ & \le \, C \,
(2^j \, r)^{-n} \, \int_{Q}  |v(t)| \, dt \le \, C \, 2^{-j n} \,
r^{-\frac n 2} \, ||v||_{L^2(Q)}. \notag
\end{align}
Thus, for all $j=0,1,2, \dots$,
\begin{equation}\label{E:2.20e}
|| \nabla ( \eta_j \, \partial_i \,  \Delta^{-1} \, v) ||_{L^2(2^{j+1} Q)} \le
C \, 2^{-j \frac n 2} \, ||v||_{L^2(Q)}.
\end{equation} 
which proves (\ref{E:2.20b}).

The proof of (\ref{E:2.20c})
 for $j=2, 3, \ldots$, provided $\int_Q v(x) \, dx =0$,
is similar to that 
of (\ref{E:2.20}).   Using estimates (\ref{E:esta}) and (\ref{E:estb}), we
deduce that, for $x\in 2^{j+1} Q \setminus 2^{j-1} Q$,
\begin{align}
 & | \nabla
( \eta_j \, \partial_i \, \, \Delta^{-1} \, v)(x)|  \le
 | \nabla \eta_j (x)| \, | \partial_i \, \Delta^{-1} \, v (x)|
+ |\eta_j(x)| \, |\nabla  \partial_i \, \Delta^{-1} \, v (x)| 
 \notag \\ & \le \, C \, (2^j r)^{-1} \,
 \int_{Q} |K^i (x-t) - K^i(x)| \, |v (t)| \, dt  \notag \\ &
+ C \, \sum_{m=1}^n \, 
 \int_{Q} | K^{i, \, m} (x-t) - K^{i, \, m} (x)| \, |v (t)| \, dt
  \le \, C \, 2^{-j(n+1)} \, |Q|^{-1} \, \int_{Q}  |v(t)| \, dt. \notag
\end{align}
This yields
$$
|| \nabla ( \eta_j \, \partial_i \, \, \Delta^{-1} \, v)||_{L^2(2^{j+1} Q)}
\le
C \, 2^{-j(1 + \frac n 2)} \, |Q|^{-\frac 1 2} \, ||v||_{L^2(Q)}.\notag
$$

The proof of (\ref{E:2.20d})
 for $j=2, 3, \ldots$ is very similar to that of (\ref{E:2.20c}),  and
is omitted here. \end{proof}

For $N>0$, define a smooth cut-off function $\psi_N(x) = \psi(\frac {x} N)$, 
where 
\begin{equation}\label{E:cutoff} 
\psi \in C^\infty_0 (\R^n); \quad \psi (x) = 1 \, \, {\rm if} \, \,  
|x|\le \frac 1 2, \quad \psi(x)=0 \, \, {\rm if} \, \,  |x|>1.  
\end{equation}

\begin{lemma} \label{Lemma 3.3} Suppose 
$\vec b \in D'(\R^n)^n$, 
$n \ge 2$. Suppose that {\rm (\ref{E:3.4q})} holds. Then 
\begin{equation}\label{E:hodge} 
\vec b = \nabla f + {\rm Div} \, F \quad {\rm in}  \, \, D'(\R^n)^n,
\end{equation}
where 
\begin{equation}\label{E:bmo} 
f=\Delta^{-1} {\rm div} \, \vec b \in {\rm BMO}(\R^n), \qquad 
F=\Delta^{-1} {\rm curl} \, \vec b \in {\rm BMO}(\R^n)^{n \times n}.  
\end{equation}
Here $f$ and $F$ are defined (up to a constant) by, respectively, 
\begin{align}\label{E:defa} 
& f  = \lim_{N \to + \infty} \, f_N, \qquad 
 f_N=\Delta^{-1} {\rm div} \, (\psi_N \vec b), \\ 
& F  = \lim_{N \to + \infty} \, F_N,  \qquad 
F_N = \Delta^{-1} {\rm curl} \, (\psi_N \vec b), \label{E:defb} 
\end{align}
 in the sense of the convergence in the 
weak-$*$ topology of ${\rm BMO}(\R^n)$. The limits above do not depend 
on the choice of $\psi_N$. 

Furthermore,  
\begin{align}
& \nabla f  = \lim_{N \to + \infty} \nabla f_N, \quad 
{\rm Div} \, F = \lim_{N \to + \infty} {\rm Div} \, F_N \quad {\rm in} 
\quad D'(\R^n)^n, \label{E:conv} \\ 
 & {\rm curl}  \, (\nabla f)   = 0, \quad
{\rm div} \, ({\rm Div} \, F) = 0, \quad   
\Delta f  = {\rm div} \, \vec b,  \quad   
 \Delta F = {\rm curl} \, \vec b.\label{E:lapl}  
\end{align}
 \end{lemma}

\begin{proof}  By 
Proposition~\ref{Proposition 3.2},  (\ref{E:3.3})  implies 
 (\ref{E:3.4q}). It follows that the latter 
inequality holds with $\psi_N \, \vec b$ in place of $\vec b$, i.e.,
for every cube $Q$, 
\begin{equation}\label{E:3.6a}
|| \psi_N \, \vec b||_{L^{-1, \, 2}(Q)} \le C \, |Q|^{\frac 1 2},  
\end{equation}
where $C$ does not depend on $Q$ and $N$. This is a consequence of 
the inequality  
$
||(\nabla \psi_N) \, v||_{L^2(\R^n)}  \le C(n) \, 
||\nabla  v||_{L^2(\R^n)},
$ for $v \in C^\infty_0(\R^n)$, 
which follows from Poincar\'e's inequality.  

We observe that  $f_N$ and $F_N$ given respectively by  (\ref{E:defa}) 
and  (\ref{E:defb}) 
are well-defined in terms of distributions. Moreover, 
by (\ref{E:3.6a}),  $\psi_N \, \vec b \in L^{-1, \, 2}(\R^n)$, 
and hence $f_N \in L^{2}(\R^n)$, $F_N \in L^{2}(\R^n)^{n\times n}$. 

Our next step is to show that, for $i, m =1, 2, \ldots, n$, 
\begin{equation}\label{E:3.7a}
||  \partial_i \, \partial_m \, \Delta^{-1} \,
 (\psi_N \, \vec b) ||_{L^{-1, \, 2} (Q)} \le C \, |Q|^{\frac 1 2},  
\end{equation}
where $C$ does not depend on $Q$ and $N$.

Notice that  $\partial_i \, \partial_m \, 
\Delta^{-1} \, (\psi_N \, \vec b)$ is well-defined in $D'(\R^n)^n$. 
Let $\{\eta_j\}_{j=0}^\infty$ be the 
partition of unity  (\ref{E:res1})$-$(\ref{E:res2}) 
associated with a cube $Q$ in $\R^n$. Then 
$$
\langle \partial_i \, \partial_m \, \Delta^{-1} \, (\psi_N \, \vec b), \,
\vec v \rangle=
\langle \vec b, \,  \psi_N \, \Delta^{-1} \, 
\partial_i \, \partial_m \vec v \rangle = 
\sum_{j=0}^\infty \, \langle \psi_N \, \vec b, \,
 \eta_j \,  \partial_i \, \partial_m \, \Delta^{-1}  \vec v \rangle,
$$
for every $\vec v  \in C^\infty_0(Q)^n$, where the sum on the 
right contains only a finite number of nonzero terms. 

 Then by  (\ref{E:3.6a})
and  statement (i) of Proposition~\ref{Proposition 2.2}, 
\begin{align}
\left \vert \langle \psi_N \, \vec b, \,   \partial_i \, \partial_m \, 
\Delta^{-1} \, \vec v \rangle \right \vert & \le \sum_{j=0}^{\infty} \,
\left \vert \langle \psi_N \, \vec b, \, \eta_j \, \partial_i \, \partial_m \,
\Delta^{-1} \,
\vec v \rangle \right \vert
\notag\\& \le c \, \sum_{j=0}^{\infty} \, 2^{j \frac n 2} \, |Q|^{\frac 1 2}
 \, || \nabla
( \eta_j  \, \partial_i \, \partial_m \, \Delta^{-1} \, \vec v)||_{L^2(2^{j+1} Q)}
\notag\\& \le C \, |Q|^{\frac 1 2} \,  ||\nabla v||_{L^2(Q)},
\notag 
\end{align}
 i.e.,  (\ref{E:3.7a}) holds. In particular, 
$$
 ||\nabla \,  f_N||_{L^{-1, \, 2} (Q)} 
\le C \,  |Q|^{\frac 1 2},  \qquad 
||\mathbf D (F_N) ||_{L^{-1, \, 2} (Q)}
\le C \,  |Q|^{\frac 1 2}.
$$
This gives:
\begin{align}
||f_N - m_Q(f_N)||^2_{L^2(Q)} & \le c \,
 || \nabla f_N||^2_{L^{-1, \, 2}(Q)} \le   C \,  |Q|, \notag \\ 
||F_N - m_Q(F_N)||^2_{L^2(Q)} & \le c \,
 || \mathbf D (F_N)||^2_{L^{-1, \, 2}(Q)} \le   C \,  |Q|,\notag
\end{align}
 where $C$ does not depend on $Q$ and $N$. Hence, 
$$\sup_{N} \, ||f_N||_{{\rm BMO}(\R^n)} < \infty, \qquad 
\sup_N \, ||F_N||_{{\rm BMO}(\R^n)^{n \times n}} < \infty.
$$

We now show that both $\{f_N\}$ and 
$\{F_N\}$ converge in the weak-$*$ topology of ${\rm BMO}$
(considered as the
dual of $\mathcal{H}^1$; see \cite{St}) respectively to 
$f \in {\rm BMO} (\R^n)$, and 
$F \in {\rm BMO} (\R^n)^{n\times n}$
(defined up to an additive constant).  We will then 
deduce that $\Delta f = {\rm div} \,  \vec b$ and  
$\Delta F = {\rm curl} \,  \vec b$ in the distributional sense, 
and   set 
$$
f = \Delta^{-1} {\rm div} \,  \vec b, \qquad 
F = \Delta^{-1} {\rm curl} \,  \vec b. 
$$

Let us prove the weak-$*$   convergence for the sequence $\{f_N\}$ in 
${\rm BMO} (\R^n)$. 
For $\{F_N\}$, the 
argument  is quite similar. Since $\{f_N\}$ is uniformly bounded 
in the BMO-norm, it is 
enough to verify that it forms a Cauchy sequence 
in the weak-$*$ topology of ${\rm BMO}$ on a dense family of
$C^\infty_0$-functions in $\mathcal{H}^1(\R^n)$.

Suppose that $v\in C^\infty_0(\R^n)$ is supported
on a cube $Q$, and 
$\int_Q v(x) \, dx =0$. Then using the same partition of unity 
$\{\eta_j\}$ one can easily check that 
$$\left \vert \int_{\R^n} (f_N - f_M) \, \bar v \, dx \right \vert
\le \sum_{j \ge N_0} \, \left \vert \langle
(\psi_N - \psi_M) \, \vec b, \,  \eta_j \,
\nabla \, \Delta^{-1} \, v \rangle\right \vert,
$$
where $N_0 \to +\infty$ as $M, N \to +\infty$.  We deduce
from (\ref{E:3.6a}):
$$\left \vert \langle (\psi_N - \psi_M) \, \vec b, \,  \eta_j \,
\nabla \, \Delta^{-1} \, v \rangle\right \vert
\le c \, |2^j \, Q|^{\frac 1 2}
\, ||\nabla  ( \eta_j \,
\nabla  \, \Delta^{-1} \, v) ||_{L^2 (2^j Q)}.$$
By statement (iii) of 
 Proposition~\ref{Proposition 2.2}, 
\begin{equation}\label{E:3.9a}
 ||\nabla  ( \eta_j \, \nabla \, \Delta^{-1} \, v)||_{L^2 (2^j Q)}  
\le c \, 2^{-j (1+\frac n 2)} |Q|^{-\frac 1 2} \, ||v||_{L^1(Q)},
\quad j \ge N_0,
\end{equation}
where $c$ does not depend on $j$, $Q$, and $v$.
Thus,
$$\left \vert \langle (\psi_N - \psi_M) \, \vec b, \,  \eta_j \, 
\nabla \, \Delta^{-1} \, v \rangle\right \vert
\le c \, 2^{-j} \, 
||v||_{L^1(Q)}, \quad j \ge N_0,
$$
and consequently,
$$
\sum_{j \ge N_0} \, \left \vert \langle (\psi_N - \psi_M) \, \vec b, 
\,  \eta_j \, 
\nabla \, \Delta^{-1} \, v \rangle\right \vert
\le c \, ||v||_{L^1(Q)} \, \sum_{j \ge N_0} \,  2^{-j}, \quad j \ge N_0.
$$

Using the preceding estimates and letting $M, \, N \to +\infty$
so that $N_0 \to +\infty$, we see that $\{f_N\}$ is 
a Cauchy sequence  in the weak-$*$ topology of ${\rm BMO}$.
In particular,
\begin{equation}\label{E:3.10a}
 \lim_{N \to + \infty} \int_{\R^n} f_N \, \bar v  \, dx  =
 \int_{\R^n} f \, \bar v \, dx, \qquad   v \in
C^\infty_0 (\R^n), \, \,
\int_{\R^n} v \, dx =0,
\end{equation}
where $f \in {\rm BMO}(\R^n)$.

To show that the limit in (\ref{E:3.10a})  does not depend on the choice of 
the cut-off functions $\psi_N$, 
and for future reference, we now demonstrate 
 that, for every $v \in C^\infty_0(\R^n)$ 
 supported on a cube $Q$ such that  
$\int_{Q} v \, dx =0$, it follows 
 \begin{equation}\label{E:3.11a}
 \int_{\R^n} f \, \bar v  \, dx  = - \sum_{j= 0}^\infty \,  \langle
 \vec b, \,  \eta_j \, \nabla (\Delta^{-1} \, v) \rangle.
\end{equation}

Notice that, by (\ref{E:3.4q}) and statement (iii) of 
Proposition~\ref{Proposition 2.2}, 
\begin{align}
\sum_{j \ge M} \,  \left \vert \langle
 \vec b, \,  \eta_j \, \nabla \Delta^{-1} \, v \rangle\right \vert 
& \le C \, \sum_{j \ge M} \, |2^j Q|^{\frac 1 2} \, 
||\nabla (\eta_j \, \nabla \Delta^{-1} \, v)||_{L^2 (2^j Q)} \notag 
\\ & \le C \, ||v||_{L^1(Q)} \, \sum_{j \ge M} \, 2^{-j},
\notag
\end{align}
for every $M \ge 1$. Moreover, by (\ref{E:3.6a}), a similar estimate holds 
with $\psi_N \vec b$ in place of $\vec b$, and $C$ which does not 
depend on $M$ and $N$. 

Clearly, (\ref{E:3.11a}) holds with $\psi_N \vec b$ 
in place of $\vec b$, and,  for $N$ large, 
$$
\sum_{0 \le j \le M} \,  \langle
  \vec b, \,  \eta_j \, \nabla \Delta^{-1} \, v  \rangle  = 
\sum_{0 \le j \le M} \, \langle
 \psi_N   \vec b, \,  \eta_j \, \nabla \Delta^{-1} \, v \rangle.
$$
By picking $M$ and $N$ 
large enough, and taking into account the above estimates together 
with (\ref{E:3.10a}), we arrive at (\ref{E:3.11a}). 

We observe that  (\ref{E:3.11a}) with 
${\rm div} \, \vec v $ in place of $v$ yields
 \begin{equation}\label{E:repr}
  \langle \nabla f, \, \vec v\rangle = 
- \int_{\R^n} f \, \, \overline{ {\rm div} \,  \vec v}  \, dx  = 
\sum_{j= 0}^\infty \,  \langle
 \vec b, \,  \eta_j \, \nabla (\Delta^{-1} {\rm div} \, \vec  v) \rangle,
\end{equation}
for every $\vec v \in C^\infty_0(\R^n)^n$ supported on $Q$. 
Hence, $\nabla f \in 
D'(\R^n)^n$, and 
$$ 
\nabla  f = \lim_{N \to + \infty} \,
\nabla f_N  \quad {\rm  in} \, \,
D'(\R^n), \qquad {\rm curl}  \, (\nabla f)   = 0 \quad {\rm in} \, \,  
 D'(\R^n)^{n \times n}. 
$$
 Moreover, for every $v \in C^\infty_0 (\R^n)$, 
$$
 \langle \Delta f, \, v \rangle   
= \lim_{N \to + \infty} \, \langle f_N, \, \Delta v \rangle 
=-  \lim_{N \to + \infty} \, \langle \psi_N \vec b, \, \nabla v \rangle 
= - \langle \vec b, \, \nabla v \rangle, 
$$
which gives $\Delta f = {\rm div} \, \vec b$ in $D'(\R^n)$.

In a completely analogous fashion, one verifies that 
$F_N \to F$ in the weak-$*$ topology of ${\rm BMO}$, 
$$ 
\lim_{N \to + \infty} \,
{\rm Div} \, F_N =  {\rm Div} \, F  \quad {\rm  in} \, \,
D'(\R^n)^n,
$$
and $\Delta F = {\rm curl} \, \vec b$ in $D'(\R^n)^{n \times n}$, 
${\rm div} \, ({\rm Div} \, F)=0$. Moreover, 
$F$ is a skew-symmetric matrix field since $F_N$ is skew-symmetric  
for every $N$.

We are now in a position to establish decomposition (\ref{E:hodge}) 
for vector fields which obey (\ref{E:3.3}).  
Let us set 
$\vec c= \nabla \, f$ and $\vec d = {\rm Div} \, F.$
Using a standard decomposition for $\vec v \in C^\infty_0(\R^n)^n$, 
\begin{equation}\label{E:hodgec}
\vec v = \nabla ( \Delta^{-1} {\rm div} \, \vec v) +
{\rm Div} \, ( \Delta^{-1} {\rm curl} \, \vec v),
\end{equation}
we deduce:
\begin{align} 
 \langle \nabla f_N, \, \vec v \rangle  
& = - \langle f_N, \, {\rm div} \, \vec v \rangle = 
\langle \psi_N \, \vec b, \, \nabla ( \Delta^{-1} {\rm div} \, 
\vec v) \rangle   \notag \\ & =  \langle \psi_N \, \vec b, \,
 \vec v \rangle - 
\langle \psi_N \, \vec b, \, {\rm Div} \,  ( \Delta^{-1} {\rm curl} \, 
\vec v) \rangle.\notag
\end{align}
Hence, 
\begin{align} 
\langle \vec c, \, \vec v \rangle   & = 
\lim_{N \to +\infty}
 \langle \, \nabla f_N,  \vec v\rangle 
 = \lim_{N \to +\infty}  \langle \psi_N \, \vec b, \, \vec v \rangle  
- \lim_{N \to +\infty} \langle \, \psi_N \, \vec b, \, {\rm Div} \, 
 (\Delta^{-1} {\rm curl} \, \vec v) \rangle\notag \\ 
& = \langle \vec b, \, \vec v \rangle - \lim_{N \to +\infty}
 \langle \, {\rm Div} \, F_N,  \vec v\rangle  
= \langle \vec b, \, \vec v \rangle - \langle \vec d, \, \vec v \rangle.
\notag
\end{align}
This completes the proof of Lemma~\ref{Lemma 3.3}. 
\end{proof}

\begin{corollary} \label{Corollary 3.4} 
Denote by $P$ and $Q$ respectively the operators 
\begin{equation}\label{E:proj} 
P = \nabla ( \Delta^{-1} {\rm div} ), \qquad Q = 
{\rm Div} \,  ( \Delta^{-1} {\rm curl} )
\end{equation}
defined on the class of vector fields $\vec b$ which 
obey {\rm (\ref{E:3.3})}. 
Then $P$ and $Q$ are bounded complementary projections, that is,  
 both $P \vec b$ and 
$Q \vec b$ satisfy {\rm (\ref{E:3.3})},  $P (P \vec b) = P \vec b$, 
$Q (Q \vec b) = Q \vec b$, and $P \vec b + Q \vec b = \vec b$. 
 \end{corollary}

\begin{proof} Suppose $\vec b \in D'(\R^n)$, and  {\rm (\ref{E:3.3})} 
holds. Let $\vec c = P \vec b$ and $\vec d = Q \vec b$. 
By Lemma~\ref{Lemma 3.3}, $ \vec c + \vec d = \vec b$. Moreover,  
${\rm curl} \, \vec b = {\rm curl} \, \vec d$, and $\vec d = {\rm Div} \, 
F$, where 
$F=\Delta^{-1} {\rm curl} \, \vec d \in {\rm BMO}(\R^n)^{n \times n}$.  
 
Then,  for every $u$, $v \in C^\infty_0(\R^n)$,
$$\left \vert \langle \vec d, \, \bar u \nabla v - v \nabla \bar u 
\rangle \right\vert = \left \vert 
{\rm trace} \, \langle F, \, \mathbf D[\bar u \nabla v - v \nabla \bar u] 
\rangle \right\vert
\le C \, ||\nabla u||_{L^2(\R^n)} \, ||\nabla v||_{L^2(\R^n)},
$$
 by the div-curl lemma \cite{CLMS}. (If $n=2$, this is equivalent to 
the Jacobian estimate
in $\mathcal H^1(\R^2)$; for $n\ge 3$, it follows by 
the commutator 
estimates involving Riesz transforms. See \cite{T}, Sec. 3.8, and 
the proof of Theorem~\ref{Theorem 3.3} below.)  
Thus, {\rm (\ref{E:3.3})} holds with $\vec d$, and hence $\vec c$,  
in place of $\vec b$.  

It remains to verify that 
 $P (P \vec b) = P \vec b$. By the preceding estimate and 
Proposition~\ref{Proposition 3.2} applied to $\vec d$,  it follows 
\begin{equation}\label{E:cur}
||\vec d||_{L^{-1, \, 2}(Q)} \le C \, |Q|^{\frac 1 2},
\end{equation}
for every cube $Q$. Then, obviously,
\begin{equation}\label{E:cur1}
||\nabla \psi_N \cdot\vec d||_{L^{-1, \, 2}(Q)} 
\le C \, N^{-1} \,  |Q|^{\frac 1 2}.
\end{equation}
where $C$ does not depend on $Q$ and $N>0$. 

From this we will deduce 
\begin{equation}\label{E:zero}
\lim_{N \to +\infty}  \, \langle \nabla \psi_N \cdot \vec d, \,
\Delta^{-1} \,  {\rm div} \,  \vec v \rangle  =0, \qquad  
\vec v \in C^\infty_0(\R^n)^n.
\end{equation}
Observe that $\nabla \psi_N (x)=0$ 
unless $\frac N 2 \le |x| \le N$, 
and thus, for $\vec v \in C^\infty_0(Q)^n$, 
$$
\langle  \nabla \psi_N \cdot \vec d,
\,  \Delta^{-1} {\rm div} \, \vec v \rangle = \sum_{M_1 \le j \le M_2} \, 
\langle  \nabla \psi_N \cdot \vec d, \, \eta_j 
\,  \Delta^{-1} {\rm div} \, \vec v \rangle,
$$
where $M_1$, $M_2 \to +\infty$ as $N \to + \infty$. 
Recall that $\eta_j$ is supported on $2^{j+1}Q \setminus 
2^{j-1} Q$ for $j \ge 1$. Hence, $\nabla \psi_N \, \eta_j$ is supported on
 $ \{2^{j+1}Q \setminus 2^{j-1} Q \}\cap 
\{\frac N 2 \le |x| \le N\}$. We may assume without loss of generality 
that $|a_Q|< 2^{j} \ell(Q)$ for $N$ large, where $a_Q$ denotes the center of $Q$. 
Then clearly, $c_1(n) \frac N {\ell(Q)} 
\le 2^j \le  c_2(n) \frac {N} {\ell(Q)}$. In other words, for a fixed $Q$, 
the sum above contains a bounded number of terms which does not depend 
on $N$.  

Thus, by (\ref{E:cur1}),  
\begin{align}
& \left \vert \langle  \nabla \psi_N \cdot \vec d,
\,  \Delta^{-1} {\rm div} \, \vec v \rangle \right \vert  \le 
\sum_{M_1 \le j \le M_2} \, 
\vert \langle  \nabla \psi_N \cdot \vec d, \, \eta_j 
\,  \Delta^{-1} {\rm div} \, \vec v \rangle \vert \notag \\ & \le C \, N^{-1} \, 
|2^{j+1} Q|^{\frac 1 2} \, \sum_{M_1 \le j \le M_2} \, 
||\nabla (\eta_j \Delta^{-1} {\rm div} \, \vec v)||_{L^2(2^{j+1} Q)}. 
\notag
\end{align} 
By statement (ii) of Proposition~\ref{Proposition 2.2}, 
$$
||\nabla (\eta_j \Delta^{-1} {\rm div} \, \vec v)||_{L^2(2^{j+1} Q)} \le C \, 
2^{-j \frac n 2} \,  ||\vec v||_{L^2(Q)}.
$$
Combining the preceding estimates we obtain 
$$
 \left \vert \langle  \nabla \psi_N \cdot \vec d,
\,  \Delta^{-1} {\rm div} \, \vec v \rangle \right \vert \le C \, N^{-1} \, 
|Q|^{\frac 1 2} \, 
 ||\vec v||_{L^2(Q)},
$$
and hence (\ref{E:zero}) holds.

By Lemma~\ref{Lemma 3.3}, ${\rm div} \, \vec d = 0$, and so integration by 
parts yields
\begin{align}
\langle \,  \psi_N \, \vec d, \, \nabla ( \Delta^{-1}
 {\rm div} \, \vec v) \rangle & =
 - \langle  \, \nabla \psi_N \cdot \vec d,  \, 
\Delta^{-1} {\rm div} \, \vec v \rangle - 
\langle \, \psi_N \, {\rm div} \, \vec d,  \, \Delta^{-1} 
{\rm div} \, \vec v \rangle, \notag \\ & 
= - \langle  \, \nabla \psi_N \cdot \vec d, 
\, \Delta^{-1} {\rm div} \, \vec v\rangle.
\notag
\end{align}

Thus,
$$
\langle P \vec d, \, \vec v\rangle = 
\lim_{N\to +\infty} \, \langle  \psi_N \, \vec d,
\,  \nabla \, (\Delta^{-1} {\rm div} \, \vec v) \rangle = 0,
$$
i.e.,  $P \vec d = 0$. Consequently, 
$P (P \vec b) = P (\vec b -\vec d) = P \vec b$, i.e., $P$, and hence 
$Q$, is a projection. 
\end{proof}

\begin{lemma} \label{Lemma 3.5} Suppose $\vec b \in D'(\R^n)^n$ 
and $q \in D'(\R^n)$, $n \ge 2$.  
Suppose  $\vec b \cdot \nabla + q$ is form bounded on 
$L^{1, \, 2}(\R^n)$, i.e., both {\rm (\ref{E:3.2})} and 
{\rm (\ref{E:3.3})} hold. Then  $q = {\rm div} \, \vec b = 0$ if $n=2$.
 In the case $n \ge 3$, 
\begin{equation}\label{E:bmo1} 
\Delta^{-1} {\rm div} \, \vec b \in {\rm BMO} (\R^n), \qquad 
\Delta^{-1} q \in {\rm BMO} (\R^n), 
\end{equation} 
where $\Delta^{-1} {\rm div} \, \vec b$ and $\Delta^{-1} q$ 
are  defined (up to a constant) by, respectively, 
\begin{equation}\label{E:lim} 
\Delta^{-1} {\rm div} \, \vec b  = 
\lim_{N \to \infty} \, 
\Delta^{-1} (\psi_N \, {\rm div} \, \vec b), 
\qquad  \Delta^{-1} q = \lim_{N \to \infty} \, 
\Delta^{-1} (\psi_N \, q),
\end{equation}
in terms of the convergence in the weak-$*$ topology of ${\rm BMO} (\R^n)$.

Furthermore, 
for every cube $Q$ in $\R^n$, 
\begin{equation}\label{E:3.5c} 
\int_Q \left ( |\nabla (\Delta^{-1} {\rm div} \, \vec b)|^2 + 
|\nabla (\Delta^{-1} q)|^2 \right ) 
\, dx \le C \, |Q|^{1- \frac 2 n}, 
\end{equation}
 where $C$ does not depend on $Q$. 
\end{lemma}

\begin{remark}\label{Remark 4} We have already defined 
$\Delta^{-1} {\rm div} \, \vec b$ in Lemma~\ref{Lemma 3.3} by  
$\Delta^{-1} {\rm div} \, \vec b = \lim_{N \to \infty} \, 
\Delta^{-1} {\rm div} \, (\psi_N \, \vec b)$. 
However, as we will show below,  this definition is consistent 
with (\ref{E:lim}) under the assumptions of Lemma~\ref{Lemma 3.5}. 
\end{remark} 

\begin{proof} By (\ref{E:3.2}), 
$$
\vert \langle q - \tfrac 1 2 \, {\rm div} \, \vec b, \, \bar u \, v 
\rangle \vert \le C \,
||\nabla u||_{L^2(\R^n)} \, ||\nabla v||_{L^2(\R^n)},
\qquad u, v \in C^\infty_0(\R^n).
$$
Letting $v \in C^\infty_0(Q)$ and $u \in C^\infty_0(2Q)$, $u=1$ on $Q$, 
in the preceding inequality, and taking the infimum over all such $u$,
as in the proof of Proposition~\ref{Proposition 3.2},
we arrive at the estimate
\begin{equation}\label{E:q}
\vert \langle q - \tfrac 1 2 \, {\rm div} \, \vec b, 
\, v \rangle \vert \le C \, {\rm cap} \, (Q)^{\frac 1 2} 
 \, ||\nabla v||_{L^2(Q)},
\qquad v \in C^\infty_0(Q).
\end{equation}

In two dimensions, ${\rm div} \, \vec b=0$ by 
Proposition~\ref{Proposition 3.2}, and 
  ${\rm cap} \, (Q) =0$ by (\ref{E:cube}). Hence, we see 
from (\ref{E:q}) that $\langle q, \, v \rangle = 0$ for every 
$v \in C^\infty_0(\R^n)$, i.e., $q=0$.

Let us now consider the case $n \ge 3$. By  (\ref{E:q}),  
$$
\vert \langle q - \tfrac 1 2 \, {\rm div} \, \vec b,  v \rangle 
\vert \le C \,
|Q|^{\frac 1 2 -\frac 1 n} \, ||\nabla v||_{L^2(\R^n)}, \qquad
v \in C^\infty_0(Q).
$$
Notice that by Proposition~\ref{Proposition 3.2},
\begin{equation}\label{E:loc}
\vert \langle {\rm div} \, \vec b,  v \rangle \vert \le C \,
|Q|^{\frac 1 2 -\frac 1 n} \, ||\nabla v||_{L^2(\R^n)}, \qquad
v \in C^\infty_0(Q).
\end{equation}
Combining the preceding estimates, we obtain 
\begin{equation}\label{E:q1}
\vert \langle q,  v \rangle 
\vert \le C \,
|Q|^{\frac 1 2 -\frac 1 n} \, ||\nabla v||_{L^2(\R^n)}, \qquad
v \in C^\infty_0(Q).
\end{equation}

Thus, 
\begin{equation}\label{E:q1a}
|| {\rm div} \, \vec b||_{L^{-1, \, 2} (Q)} + 
|| q||_{L^{-1, \, 2} (Q)}\le C \, |Q|^{\frac  1 2 - \frac 1 n}.  
\end{equation}
 This obviously implies
\begin{equation}\label{E:q1b} 
||\psi_N \, {\rm div} \, \vec b||_{L^{-1, \, 2} (Q)} 
+ || \psi_N \, q||_{L^{-1, \, 2} (Q)} 
\le C \, |Q|^{\frac  1 2 - \frac 1 n}, 
\end{equation}
where
$C$ does not depend on $Q$ and $N$. 

We now set
$$
\tilde f_N = \Delta^{-1} (\psi_N \, {\rm div} \, \vec b), 
\qquad g_N= \Delta^{-1} (\psi_N \, q),
$$
which are well-defined in $D'(\R^n)$. Note that $\tilde f_N$ 
differs slightly 
from $f_N = \Delta^{-1}  {\rm div} (\psi_N \, \vec b)$ 
used in Lemma~\ref{Lemma 3.3}. We will deduce from  (\ref{E:q1b}) 
\begin{equation}\label{E:q1c}
|| \nabla  \tilde f_N||_{L^{-1, \, 2} (Q)} 
+ || \nabla g_N ||_{L^{-1, \, 2} (Q)} 
\le C \, |Q|^{\frac  1 2},
\end{equation}
where $C$ does not depend on $Q$ and $N$. It is enough to estimate 
only the 
second term on the left-hand side;  the first one 
 is treated analogously simply by putting ${\rm div} \, \vec b$ in place 
of $q$.

Notice that,
for every $\vec v \in C^\infty_0(Q)^n$,
$$
\langle \nabla \Delta^{-1} (\psi_N \, q),  \vec v \rangle =
-  \langle q, \,  \psi_N \, \Delta^{-1}  {\rm div} \, \vec v \rangle =
- \sum_{j=0}^{\infty} \,
\langle \psi_N \, q,  \eta_j \, \Delta^{-1}  {\rm div} \, \vec v \rangle,
$$
where the sum on the right contains only a finite number of nonzero terms. 
Now using (\ref{E:q1b}) and statement (ii) of 
Proposition~\ref{Proposition 2.2},
we estimate
\begin{align}
& \vert \langle \nabla \Delta^{-1} (\psi_N \, q),  \vec v \rangle \vert
\le C \, \sum_{j=0}^{\infty} \, \vert
\langle \psi_N \, q,  \eta_j \, \Delta^{-1}  {\rm div} \, \vec v \rangle
\vert \notag \\ & \le C \, \sum_{j=0}^{\infty} \,|2^{j+1} Q|^{\frac 1 2 - \frac 1 n} \,
||\nabla (\eta_j \, \Delta^{-1}  {\rm div} \, \vec v)||_{L^2(2^j Q)} \notag \\
& \le C \,  |Q|^{\frac 1 2 - \frac 1 n} \, ||\vec v||_{L^2(Q)} \, 
\sum_{j=0}^{\infty} \, 2^{-j}
\le C \,  |Q|^{\frac 1 2} \, ||\nabla \vec v||_{L^2(Q)}.
\notag
\end{align}
This proves (\ref{E:q1c}), from which it  is immediate  that
\begin{align}
& || \tilde f_N -m_Q( \tilde f_N) ||^2_{L^2(Q)} + 
|| g_N -m_Q(g_N) ||^2_{L^2(Q)} \notag \\ \le C \, & 
( ||\nabla \tilde f_N||^2_{L^{-1, \, 2}(Q)} +  
||\nabla g_N||^2_{L^{-1, \, 2}(Q)} ) \, \le C \, |Q|.
\notag
\end{align}
Thus, 
$$\sup_N \, (|| \tilde f_N||_{{\rm BMO}(\R^n)} +  
|| g_N||_{{\rm BMO}(\R^n)} ) < +\infty.
$$

We can now follow the argument of Lemma~\ref{Lemma 3.3}  to
demonstrate  that 
$$\tilde f_N \to \tilde f = \Delta^{-1} {\rm div} \, \vec b, \qquad 
g_N \to g = \Delta^{-1} q,$$
 in the sense of the weak-$*$ topology of ${\rm BMO} (\R^n)$. 
Note that we also have to verify $\tilde f = f$ where 
$f = \lim_{N \to \infty} \, 
f_N$. 

Let us indicate some 
changes that are needed here. We are now utilizing 
estimates (\ref{E:q1b}) in place of (\ref{E:3.6a}). 
By using  statement (iv), rather than statement (iii), 
of Proposition~\ref{Proposition 2.2}, 
one deduces that $\{g_N\}$ is a Cauchy sequence. Hence, $g$
is defined (up to a constant) by 
$$
g=\Delta^{-1} q = \lim_{N \to \infty} \, \Delta^{-1} (\psi_N \, q)
\in {\rm BMO} (\R^n),  
$$
in the sense of the weak-$*$ ${\rm BMO}$ convergence. Moreover,
using (\ref{E:q1b}) together with statement (iv) 
of Proposition~\ref{Proposition 2.2}, we obtain, exactly as in
the proof of (\ref{E:3.11a}), that 
 for any
$v \in C^\infty_0(Q)$ such that $\int_Q v \, dx =0$,  
$$\langle \Delta^{-1} q, \, v\rangle
=  \int_{\R^n}  g \, \,  \bar v \, dx 
=  \sum_{j=0}^{\infty} \,
\langle  q, \,  \eta_j \, \Delta^{-1}  {\rm div} \,  v \rangle.
$$

From the above equations we see that  
$$
\nabla g = \lim_{N \to \infty} \, \nabla \Delta^{-1} (\psi_N \, q)
\quad {\rm in} \, \, D'(\R^n),
$$
and, for every $\vec v \in C^\infty_0(Q)^n$,
$$\langle \nabla g, \, v \rangle = 
- \int_{\R^n}  g \, \,  \overline{ {\rm div} \, \vec v} \, dx
= - \sum_{j=0}^{\infty} \,
\langle  q,  \eta_j \, \Delta^{-1}  {\rm div} \, \vec v \rangle.
$$
Hence, ${\rm div} \, \nabla g = q$ in $D'(\R^n)$.

Obviously, analogous statements hold with ${\rm div} \, \vec b$ and $\tilde f$ 
in place of $q$ and $g$, respectively. It remains only to 
justify Remark~\ref{Remark 4} above. We show that 
$\tilde f =f$ in the ${\rm BMO}$ sense, i.e., 
$$ 
\int_{\R^n}   \tilde f \, \bar v \, dx = \int_{\R^n}   f \, \bar v \, dx,
$$
for  every $v \in C^\infty_0(\R^n)$ such that $\int_{\R^n} v \, dx = 0$. 
Notice that  
$$\langle \tilde f_N, \, v \rangle = \langle f_N, \, v \rangle  - 
\langle \nabla \psi_N \cdot \vec b, \, \Delta^{-1} \, v \rangle.
$$
It suffices to check 
\begin{equation}\label{E:zerob} 
\lim_{N \to \infty} \, \langle \nabla \psi_N \cdot \vec b, \, 
\Delta^{-1} \, v \rangle =0,
\end{equation}
where $v$ is supported on a cube $Q$, and $\int_Q v \, dx =0$. 
Notice that  
$$\langle \nabla \psi_N \cdot \vec b, \, \Delta^{-1} \, v \rangle = 
\sum_{M_1 \le j  \le M_2} \, \langle \nabla \psi_N \cdot \vec b, \, 
\eta_j \, \Delta^{-1} \, v \rangle, 
$$
where $M_1, \, M_2 \to \infty$ as $N \to \infty$. As 
was shown in the proof of Corollary~\ref{Corollary 3.4}, for a fixed $Q$, 
the sum above contains a uniformly bounded number of nonzero 
terms. 

We recall that by Proposition~\ref{Proposition 3.2}, 
$||\vec b||_{L^{-1, \, 2}(Q)} \le C \, |Q|^{\frac 1 2}$, and hence 
$$
||\nabla \psi_N \cdot \vec b||_{L^{-1, \, 2}(Q)} \le C \, 
N^{-1} \, |Q|^{\frac 1 2},
$$
where $C$ does not depend on $Q$ and $N$. It follows, 
\begin{align}
\vert \langle \nabla \psi_N \cdot \vec b, \, \Delta^{-1} \, v \rangle 
\vert & \le C \, \sum_{M_1 \le j \le M_2} \, 
\langle \vert \nabla \psi_N \cdot \vec b, \, 
\eta_j \, \Delta^{-1} \, v \rangle\vert \notag \\ & \le C \, N^{-1} \, 
\sum_{M_1 \le j \le M_2} \, |2^{j+1} Q|^{\frac 1 2} \, 
||\nabla ( \eta_j \, \Delta^{-1} \, v)||_{L^2(2^{j+1} Q)}.
\notag
\end{align}  
Applying  statement (iv) of  Proposition~\ref{Proposition 2.2}, we conclude:  
$$
\vert \langle \nabla \psi_N \cdot \vec b, \, \Delta^{-1} \, v \rangle 
\vert \le C \, N^{-1} \, 
||v||_{L^1(Q)}, 
$$
which yields  (\ref{E:zerob}). Thus, $\tilde f = f$ in ${\rm BMO} (\R^n)$. 

Then, for every 
$\vec v \in C^\infty_0(Q)^n$, 
\begin{align}\label{E:reprc} 
\langle \nabla f, \, \vec v\rangle & = - \int_{\R^n} f \, 
\overline{ {\rm div} \, \vec v} \, dx = - \sum_{j\ge 0} \, 
\langle {\rm div} \,  \vec b, 
\, \eta_j \Delta^{-1} \vec {\rm div} \,  \vec v\rangle, \\
\langle \nabla g, \, \vec v\rangle & = - \int_{\R^n} g \, 
\overline{ {\rm div} \, \vec v} \, dx = - \sum_{j\ge 0} \, 
\langle q,  \, \eta_j \Delta^{-1}  {\rm div} \,  \vec v\rangle. 
\label{E:reprd}  
\end{align}
This is verified exactly as in the proof of (\ref{E:3.11a}), using 
(\ref{E:q1b}) together with statement (iv), rather than  statement (iii), 
of Proposition~\ref{Proposition 2.2}. 

We are now in a position to obtain the estimate 
\begin{equation}\label{E:3.5'}
\left \vert \langle  \nabla \Delta^{-1} \, ({\rm div} \, \vec b),  
\, \vec v 
\rangle \right \vert \le C \, |Q|^{\frac 1 2  - \frac 1 n}
  \, ||\vec v||_{L^2(Q)}, \qquad   \vec v \in C^\infty_0(Q)^n. 
\end{equation}
Indeed, by 
(\ref{E:reprc})  and 
statement (ii) of Proposition~\ref{Proposition 2.2},
\begin{align}
\vert \langle \nabla f, \, \vec v\rangle\vert & \le C \, \sum_{j\ge 0} \, 
|2^{j+1} Q| \, ||\nabla(\eta_j \Delta^{-1} \vec {\rm div} \,  
\vec v)||_{L^2(2^{j +1} Q)} 
\notag\\ & \le C \, |Q|^{\frac 1 2 -\frac 1 n} \, ||\vec v||_{L^2(Q)} \, 
\sum_{j\ge 0} \, 2^{-j}. 
 \notag
\end{align} 
Taking the supremum over all $\vec v \in C^\infty_0(Q)^n$ in (\ref{E:3.5'}), 
we obtain   
\begin{equation}\label{E:3loc}
||\nabla (\Delta^{-1}  {\rm div} \, \vec b) ||_{L^2(Q)} \le C
\, |Q|^{\frac 1 2- \frac 1 n}. 
\end{equation}
Analogously, we deduce from (\ref{E:reprd}),  
\begin{equation}\label{E:3.33loc}
||\nabla (\Delta^{-1} q)||_{L^2 (Q)}
\le C \, |Q|^{\frac 1 2- \frac 1 n}.
\end{equation}
Combining the preceding estimates, we arrive at (\ref{E:3.5c}).  
\end{proof}

We now establish the main lemma, whose proof makes use of the idea 
of the magnetic gauge invariance. 

\begin{lemma} \label{Lemma 3.6} 
Let $\vec b \in D'(\R^n)$, $n \ge 2$. Suppose that,     
for all $u, \, v \in C^\infty_0(\R^n)$, 
 \begin{equation}\label{E:3.3i}
\left \vert 
 \langle   \vec b,  \, 
 \bar u \, \nabla v  - v \nabla \bar u \rangle \right \vert 
\le C \, 
||u||_{L^{1, \, 2} (\R^n)} \, ||v||_{L^{1, \, 2}(\R^n)}. 
\end{equation} 
Then 
$\vec c = \nabla (\Delta^{-1} {\rm div} \, \vec b) 
\in L^2_{{\rm loc}}(\R^n)$, and $|\vec c|^2 \in \mathfrak{M}^{1, \, 2}_+(\R^n)$, i.e.,  
\begin{equation}\label{E:3.4i}
\int_{\R^n}  |u|^2 \, |\vec c|^2 \, dx \le C_1 \, 
||u||^2_{L^{1, \, 2} (\R^n)}, 
\end{equation} 
where $C_1$ does not depend on $u \in C^\infty_0(\R^n)$. If $n=2$, then $\vec c = 0$. \end{lemma} 

\begin{proof} Suppose that {\rm (\ref{E:3.3i})} holds. Then by 
continuity the bilinear form on the left-hand side can be 
extended to all $u, \, v \in L^{1, \, 2} (\R^n)$. 
Let $v$ be a nonnegative function such that $\lambda = \log v$ 
has the property:
$$
\nabla \lambda  = \frac {\nabla v} {v} \in L^2_{{\rm loc}} (\R^n).
$$
Moreover, we need $\lambda$ to be chosen so that 
the energy space $L^{1, \, 2} (\R^n)$ be invariant under the gauge 
transformation:  
\begin{align}
& u \to \tilde u = e^{i \, \lambda} \, u, \qquad v \to \tilde  v = 
e^{i \, \lambda} \,  v,
\label{E:g1}\\
& ||\tilde u||_{L^{1, \, 2} (\R^n)} \simeq
||u||_{L^{1, \, 2} (\R^n)}, \quad 
|| \tilde v||_{L^{1, \, 2} (\R^n)} \simeq  
||v||_{L^{1, \, 2} (\R^n)}.\label{E:g2}
\end{align}

To deduce 
{\rm (\ref{E:3.4i})}, we notice that by Theorem~\ref{Trace Theorem}, it suffices to show  
that, 
for every compact set $e\subset \R^n$,   
\begin{equation}\label{E:3.5i}
\int_{e} |\vec c|^2 \,  dx \le {\rm const} \, 
{\rm cap} \, (e). 
\end{equation} 
Without loss of generality we may assume that ${\rm cap} \, (e)>0$, 
since otherwise $|e| = 0$, and hence the preceding inequality is obvious. 
Denote by 
$\mu = \mu_e$ the equilibrium measure associated with $e$. Let 
$P \mu (x)= (-\Delta)^{-1} \mu$ denote the Newtonian potential of $\mu$. 

Suppose first that $n \ge 3$. Then 
$$P \mu(x) = c(n) \, \int_{\R^n} \frac {d \mu(y)}{|x-y|^{n-2}}, 
\qquad x \in \R^n.$$
We set 
\begin{equation}\label{E:3.7}
 \lambda = \tau \, \log (P \mu), \quad v = (P \mu)^\tau, 
\qquad  1< 2 \tau < \tfrac {n}{n-2}.
\end{equation}
We observe that $v \in L^{1, \, 2} (\R^n)$ by \cite{MV1}, Proposition 2.5. 
Clearly, $\nabla(e^{i \lambda} \, u) = 
(i \, u \, \nabla \lambda + \nabla u) \, e^{i \lambda}$.  
Consequently,  for every $u \in L^{1, \, 2}(\R^n)$, 
$$||e^{i \lambda} \, u||_{L^{1, \, 2} (\R^n)} 
\le \, ( ||u \, \nabla \lambda||_{L^{2} (\R^n)} 
+ ||u||_{L^{1, \, 2} (\R^n)}).$$
Note that $\nabla \lambda = \tau \, \frac {\nabla(P \mu)}{P \mu}$. Hence, 
by \cite{MV1}, Proposition 2.7, 
$$
||u \, \nabla \lambda||_{L^{2} (\R^n)} \le 2 \tau \, 
||u||_{L^{1, \, 2} (\R^n)}.$$
From this it follows 
\begin{equation}\label{E:3.8}
(1+ 2 \tau)^{-1} \, ||u||_{L^{1, \, 2} (\R^n)} \le 
||e^{i \lambda} \, u||_{L^{1, \, 2} (\R^n)} \le (1+ 2 \tau)
 \, ||u||_{L^{1, \, 2} (\R^n)}.
\end{equation}
Using similar estimates for $e^{-i \lambda} \, v$, we deduce  (\ref{E:g2}). 
Moreover, 
\begin{equation}\label{E:3.9bb}
||\tilde u||_{L^{1, \, 2} (\R^n)} \le (1+ 2 \tau) \, 
||u||_{L^{1, \, 2} (\R^n)}, 
\, \,  ||\tilde v||_{L^{1, \, 2} (\R^n)} \le (1+ 2 \tau) \, 
||v||_{L^{1, \, 2} (\R^n)}. 
\end{equation}
Applying (\ref{E:3.9bb}) and (\ref{E:3.3i}) with $\tilde u$ and 
$\tilde v$ in place of $u$ and $v$, we get
\begin{align}
\left \vert 
 \langle   \vec b,  \, 
 \overline{\tilde u} \, \nabla \tilde v  - 
\tilde v \nabla \overline{\tilde u} \rangle \right \vert 
& \le C \, 
||\tilde u||_{L^{1, \, 2}(\R^n)} \, 
||\tilde v||_{L^{1, \, 2}(\R^n)} \notag\\
& \le C \, (1 + 2 \tau)^2 \, 
||u||_{L^{1, \, 2}(\R^n)} \, 
||v||_{L^{1, \, 2}(\R^n)}.
\notag
\end{align}
Notice that  
$$\overline{\tilde u} \, \nabla \tilde v  - 
\tilde v \nabla \overline{\tilde u} = 
 \bar u \, \nabla v  - v \nabla \bar u - 
2 i \, \bar u \, v \, \nabla \lambda.
$$

Combining the preceding estimates, we obtain
\begin{align} 
& 2 \left \vert 
 \langle   \vec b,  \, 
 \bar  u \, v  \nabla \lambda \rangle \right \vert 
\le \left \vert \langle   \vec b,  \, \bar u \, \nabla v  - 
v \nabla \bar u \rangle \right \vert + \left \vert 
 \langle   \vec b,  \, 
 \overline{\tilde u} \, \nabla \tilde v  - 
\tilde v \nabla \overline{\tilde u} \rangle \right \vert 
\notag \\ 
& \le C \, (1 + (1+ 2 \tau)^2) \, ||u||_{L^{1, \, 2}(\R^n)} \, 
||\nabla v||_{L^{1, \, 2}(\R^n)}.
\notag
\end{align}
Observe that $v \, \nabla \lambda = \nabla v.$ 
Thus, we arrive at the inequality
\begin{equation}\label{E:3.10bb}
\left \vert 
 \langle   \vec b,  \, 
  \bar u \, \nabla v  \rangle \right \vert 
\le C \, \tfrac {1 + (1+ 2 \tau)^2} {2} \, ||u||_{L^{1, \, 2} (\R^n)} \, 
||v||_{L^{1, \, 2}(\R^n)}, 
\end{equation}
where $v = (P \mu)^\tau$. From the preceding estimate and  (\ref{E:3.3i}), 
we deduce:
\begin{equation}\label{E:3.10bc}
\left \vert 
 \langle   \vec b,  \, 
 (\bar u \, \nabla v  + v \nabla \bar u) \rangle \right \vert 
\le C \, (2 + (1+ 2 \tau)^2) \, ||u||_{L^{1, \, 2} (\R^n)} \, 
||v||_{L^{1, \, 2}(\R^n)}.
\end{equation}
This yields 
\begin{equation}\label{E:3.11bb}
\left \vert 
 \langle   ({\rm div} \, \vec b)  \, u,  \, v\rangle \right \vert 
\le C \, (2 + (1+ 2 \tau)^2) \, ||u||_{L^{1, \, 2} (\R^n)} \, 
||v||_{L^{1, \, 2} (\R^n)},
\end{equation}
where $u \in L^{1, \, 2} (\R^n)$ and $v= (P \mu)^{\tau}$. By Proposition 2.5 in \cite{MV1}, 
$$
||v||_{L^{1, \, 2} (\R^n)} = \tau (2 \tau -1)^{-\frac 1 2} 
\, {\rm cap} \, (e)^{\frac 1 2}.
$$
Hence by (\ref{E:3.11bb}),
\begin{equation}\label{E:3.11cc}
\left \vert 
 \langle   ({\rm div} \, \vec b)  \, u,  \, v\rangle \right \vert 
 \le c(\tau) \, \, C \, ||u||_{L^{1, \, 2} (\R^n)} 
\, {\rm cap} \, (e)^{\frac 1 2},
\end{equation}
where $c(\tau)$ depends only on $\tau$, and $C$ is the constant in (\ref{E:3.3i}). 

By letting $e = 2Q$ in (\ref{E:3.11cc}), where $Q$ is a cube in $\R^n$, and taking into account 
that in this case $P \mu =1$ on $Q$, and ${\rm cap} \, (2Q)
 \simeq |Q|^{\frac 1 2 - \frac 1 n},$ we see that 
\begin{equation}\label{E:3.11ee}
\left \vert 
 \langle   {\rm div} \, \vec b,  \, \bar u \rangle \right \vert 
 \le C  \, ||u||_{L^{1, \, 2} (Q)} 
\, |Q|^{\frac 1 2 - \frac 1 n},
\end{equation}
for every $u \in C^\infty_0(Q)$, $n \ge 3$. (The preceding estimate was already proved in Proposition~\ref{Proposition 3.2}.)  As in the proof of Lemma~\ref{Lemma 3.5}, estimate  (\ref{E:3.11ee}) yields  
$$
\Delta^{-1} ({\rm div} \, \vec b) = \lim_{N \to \infty} \, \tilde f_N \in {\rm BMO} (\R^n), \qquad 
\tilde f_N = \Delta^{-1} (\psi_N \, 
{\rm div} \, \vec b),$$
where the limit is understood in the sense of the  weak-$*$ convergence in ${\rm BMO}$. 
Moreover, by (\ref{E:3loc}), $\nabla \Delta^{-1}
 ({\rm div} \, \vec b) \in L^2_{{\rm loc}} (\R^n),$ and 
$$
||\nabla (\Delta^{-1} {\rm div} \, \vec b)||_{L^2(Q)} \le C \, |Q|^{\frac 1 2 - \frac 1 n},
$$
for every cube $Q$ if $n \ge 3$.

It remains to show that (\ref{E:3.11cc}), with arbitrary $u \in L^{1, \, 2}(\R^n)$ and 
$v =(P \mu)^\tau$,  
yields 
\begin{equation}\label{E:3.12bb}
\int_{e} |\nabla \Delta^{-1} \, ({\rm div} \, \vec b)|^2 \, dx \, \le C \, 
{\rm cap} \, (e),
\end{equation}
where $C$ does not depend on the compact set $e$. 

This is verified analogously  to the proof of the necessity part of Theorem 2.2 in \cite{MV1},  where $u$ and $v$ in  (\ref{E:3.11bb})
were picked exactly as above. For the sake of convenience, we outline the 
rest of the proof as follows. 

For $\vec \phi \in C^\infty_0(\R^n)^n$, we set
$$u = v^{-1} \, ({ \Delta^{-1} \, {\rm div} \, \vec \phi}), \qquad v = (P \mu)^{\tau}.$$ 
Then by Lemma 2.6 in \cite{MV1}, $u \in L^{1, \, 2} (\R^n),$ and 
$$
||u||_{L^{1, \, 2} (\R^n)} \le 
||\nabla (\Delta^{-1} \, {\rm div} \, \vec \phi)||_{L^2(\R^n, \, v^{-2} dx) } \le 
(\tau +1)^{\frac 1 2} (4 \tau +1)^{\frac 1 2} \,  ||u||_{L^{1, \, 2} (\R^n)},
$$
where 
$$||f ||_{L^2(\R^n, \, v^{-2} dx) } = \left (
\int_{\R^n} | f|^2 \, (P \mu)^{-2 \tau} \, {dx} \right)^{\frac 1 2}$$
is the weighted $L^2$-norm of $f=\nabla (\Delta^{-1} \, {\rm div} \, \vec \phi)$.
Hence, by (\ref{E:3.11bb}),
\begin{equation}\label{E:3.11dd}
\left \vert 
 \langle   {\rm div} \, \vec b,  \, \Delta^{-1} \, {\rm div} \, \vec \phi \rangle \right\vert \le c (\tau) \, C \, {\rm cap} \, (e)^{\frac 1 2} \, 
 ||\nabla (\Delta^{-1} \, {\rm div} \, \vec \phi)||^2_{L^2(\R^n, \, v^{-2} dx) }.
\end{equation}
Integrating by parts, we get
$$
\langle   {\rm div} \, \vec b,  \, \Delta^{-1} \, {\rm div} \, \vec \phi \rangle 
= - \langle  \nabla (\Delta^{-1} \, {\rm div} \, \vec b), \,  \vec \phi\rangle.$$
Thus,
$$| \langle  \nabla (\Delta^{-1} \, {\rm div} \, \vec b), \,  \vec \phi\rangle | 
\le c(\tau) \, C \, 
{\rm cap} \, (e)^{\frac 1 2} \, 
 ||\nabla (\Delta^{-1} \, {\rm div} \, \vec \phi)||^2_{L^2(\R^n, \, v^{-2} dx) }.
 $$
 
Notice that, for $\tau$ picked according to  (\ref{E:3.7}),  the weight $v^{-2}$ belongs to the 
 Muckenhoupt class $A_2(\R^n)$, and its Muckenhoupt bound depends only on $\tau$ and $n$  
  by Proposition 2.9 in \cite{MV1}. Hence, the Calderon--Zygmund operator 
  $ \nabla (\Delta^{-1} \, {\rm div})$ is bounded on $L^2(\R^n, v^{-2} dx)$, and 
 $$||\nabla (\Delta^{-1} \, {\rm div} \, \vec \phi)||_{L^2(\R^n, v^{-2} dx) } 
 \le c(\tau, n) \, ||\vec \phi||_{L^2(\R^n), \, v^{-2} dx }.$$
 Combining the preceding inequalities, we get  
 $$
 | \langle  \nabla (\Delta^{-1} \, {\rm div} \, \vec b), \,  \vec \phi\rangle | 
\le c(\tau, n) \, C \, {\rm cap} \, (e)^{\frac 1 2} \, ||\vec \phi||_{L^2(\R^n, \, v^{-2} dx) },
$$
for all $\vec \phi \in C^\infty_0(\R^n)^n$. 

Since $v^{-2}\in A_2(\R^n)$, and $C^\infty_0(\R^n)$ is dense in 
$L^2(\R^n, \, v^{-2} dx)$, we deduce from the preceding inequality:
$$\int_{\R^n} |\nabla (\Delta^{-1} \, {\rm div} \, \vec b)|^2 \, v^2 \, dx \le c(\tau, n)  \, C^2 \, 
{\rm cap} \, (e),$$
where $C$ is the constant in (\ref{E:3.3i}). 
Using the fact that $v = (P \mu)^\tau\ge 1$ $dx$-a.e. on $e$, we obtain (\ref{E:3.11dd}).

If $n=2$, then by 
Proposition~\ref{Proposition 3.2}, 
{\rm (\ref{E:3.3i})} yields   $ {\rm div} \, \vec b = 0$, 
and hence $\vec c= \nabla (\Delta^{-1} \, {\rm div} \, \vec b)=  0$.  \end{proof}

\begin{lemma}\label{Lemma 3.8}
Let $\vec b \in D(\R^n)$, $n \ge 2$. Then the inequality
\begin{equation}\label{E:comm}
\left \vert \langle  \, \vec b,  \,  \bar u \, \nabla v  -  v \, \nabla \, \bar u 
\rangle\right \vert 
\leq c \, \,  ||\nabla u||_{L^2(\R^n)} \, 
||\nabla v||_{L^2(\R^n)}, \quad   \, 
u, v \in C^\infty_0(\R^n),
\end{equation} 
 holds if and only if 
\begin{equation}\label{E:decomp}
\vec b = \vec c + {\rm Div} \, F,
\end{equation}
where $\vec c$ obeys {\rm(\ref{E:3.4i})}, and $F \in 
{\rm BMO}(\R^n)^{n \times n}$ is a skew-symmetric 
matrix field. 

Moreover, if {\rm(\ref{E:3.3i})} holds then  
 {\rm(\ref{E:decomp})} is valid 
with 
$\vec c = \nabla (\Delta^{-1} {\rm div} \, \vec b)$ 
obeying {\rm(\ref{E:3.4i})},  and 
$F = \Delta^{-1} {\rm curl} \, \vec b \in {\rm BMO}(\R^n)^{n \times n}$. 

In the case $n=2$, it follows that $\vec c =0$ in the statements above. 
\end{lemma} 

\begin{proof} Suppose first that $n \ge 3$.  To 
prove the ``if'' part, suppose that  {\rm (\ref{E:decomp})} holds, i.e., 
$\vec b = \vec c + \vec d,$ where 
$\vec d = {\rm Div} \, F$  is divergence free, 
$F \in {\rm BMO}(\R^n)^{n \times n}$, 
 and $\vec c$ satisfies   
{\rm (\ref{E:3.4i})}.  Then $\vec c \cdot \nabla$ is form bounded,  
since by Schwarz's inequality 
and {\rm (\ref{E:3.4i})}, we have:
$$
| \langle \vec c \cdot \nabla u, \, v \rangle|  
\le \, ||u \, \vec c||_{L^2(\R^n)} \, ||\nabla v||_{L^2(\R^n)} 
\le C \, \, ||\nabla u||_{L^2(\R^n)} \, ||\nabla v||_{L^2(\R^n)}. 
$$
 The preceding inequality obviously yields (\ref{E:3.3i}). 

It remains to show that  
\begin{equation}\label{E:2.14}
\left \vert \langle  \, \vec d,  \,  \bar u \, \nabla v  -  v \, \nabla \, \bar u 
\rangle\right \vert 
\leq c \, \,  ||\nabla u||_{L^2(\R^n)} ||\nabla v||_{L^2(\R^n)}, \quad   \, 
u, v \in C^\infty_0(\R^n).
\end{equation}

Let us start with the case $n = 3$. Then equivalently 
we have   $\vec d = {\rm curl} \, \vec F$, where 
$\vec F \in {\rm BMO} (\R^3)^3$. It follows 
\begin{align}
 \left \vert 
\langle  \vec d,  \,   \bar u \, \nabla v  -  v \, \nabla \, \bar u   
\rangle \right \vert  & = 2 \, \left \vert 
\langle  \, \vec F, \, \nabla \bar u \times \nabla v \rangle \right \vert  
\notag \\
& \le c \, ||\vec F||_{{\rm BMO}(R^3)^3} \, 
|| \nabla \bar u \times \nabla v||_{\mathcal{H}^1 (R^3)^3} \notag \\ 
& \le c_1 \, ||\vec F||_{{\rm BMO}(R^3)^3} \, ||\nabla u||_{L^{2} (\R^3)} \, 
||\nabla v||_{L^{2} (\R^3)}. \notag
\end{align}
The last inequality is  based on a standard compensated compactness 
argument  
using commutators with Riesz transforms  \cite{CLMS}.

Similarly, for $n \ge 3$, we have 
$\vec d = {\rm Div} \, F$, where $F = (f_{ij})_{i, j=1}^n$ is 
skew-symmetric, and hence (see Sec.~\ref{Section 2}) 
\begin{align}
\langle  \vec d,  \,   \bar u \, \nabla v  -  v \, \nabla \, \bar u   
\rangle  & = -  {\rm trace} \, \langle   F^t, \, \D 
(\bar u \, \nabla v  -  v \, \nabla \, \bar u) \rangle \notag \\ 
 & =  - \sum_{i, \, j=1}^n \langle f_{ij}, \, \partial_i \bar u  \, 
 \partial_j v -  \partial_j \bar u \, \partial_i v \rangle.
\notag
\end{align}

The inequality 
$$
\left \vert \langle f_{ij}, \, \partial_i \bar u \, 
\partial_j v - \partial_j \bar u \, 
 \partial_i v \rangle \right \vert \le c \, ||F||_{{\rm BMO}(R^n)^{n \times n}} \,
\, ||\nabla u||_{L^{2} (\R^n)} \, 
||\nabla v||_{L^{2} (\R^n)} 
$$ 
now follows again from the ${\mathcal H}^1$ estimates for commutators 
with Riesz transforms (see \cite{CLMS}).

If $n=2$, then as was indicated above 
{\rm(\ref{E:3.4i})} yields $\vec c = 0$. Hence, 
$\vec b = (\partial_2 g, \, - \partial_1 g)$, where $g \in {\rm BMO} 
(\R^2)$. Now {\rm(\ref{E:3.3i})} follows from 
the well-known $\mathcal{H}^1$ inequality for 
the Jacobian determinant \cite{CLMS}: 
$$
\left | \langle \vec b \cdot \nabla u, \, v \rangle \right| = 
\left | \langle g, \partial_2 \bar u \, 
\partial_1 v - \partial_1 \bar u \, 
\partial_2 v \rangle \right| \le C \, 
 ||u||_{L^{1, \, 2} (\R^n)} \, 
||v||_{L^{1, \, 2}(\R^n)}.
$$
This proves the ``if'' part of Lemma~\ref{Lemma 3.8}. 

To prove the converse, notice that 
 by Proposition~\ref{Proposition 3.2}
and Lemma~\ref{Lemma 3.3} it follows that one can set 
$\vec c = \nabla ( \Delta^{-1} {\rm div} \, \vec b)$ and 
$\vec d = {\rm Div} \, ( \Delta^{-1} {\rm curl} \, \vec b) \in {\rm BMO} 
(\R^n)^n$. Finally,  by Lemma~\ref{Lemma 3.6} we deduce that 
{\rm(\ref{E:3.4i})} holds. 
\end{proof} 

We are now in a position to obtain 
 the main result of this section. 

\begin{theorem} \label{Theorem 3.3} 
Let $n \ge 2$.  
Let $\mathcal{L} = \vec b \, \cdot \nabla + q$, where 
$\vec b \in D'(\R^n)$, and $q \in D'(\R^n)$. Then 
the following statements hold.

{\rm(i)} The bilinear form associated with $\mathcal{L}$ is bounded 
on $L^{1, \, 2} (\R^n) \times L^{1, \, 2} (\R^n)$ 
if and only if 
\begin{equation}\label{E:3.30}
\vec b = \vec c + {\rm Div}  \, F, \qquad q = {\rm div} \, \vec h, 
\end{equation} 
where $F\in {\rm BMO}(\R^n)^{n \times n}$ is a skew-symmetric matrix field, and $|\vec c |^2   + |\vec h|^2 \in \mathfrak{M}^{1, \, 2}_+(\R^n)$, i.e., 
\begin{equation}\label{E:3.11aa}
\int_{\R^n} ( |\vec c |^2   
 + |\vec h|^2  )  \, |u|^2 \, dx \le C \, 
||\nabla u||^2_{L^2(\R^n)}, 
\end{equation} 
for all $u \in C^\infty_0(\R^n)$.

{\rm(ii)} If  $\mathcal{L}$ is form bounded, then 
in statement {\rm (i)} 
one can set 
\begin{equation} \label{E:3.delta}
\vec c = \nabla (\Delta^{-1} {\rm div} \, \vec b), \quad
F = \Delta^{-1} {\rm curl} \, \vec b, \quad
\vec h = \nabla (\Delta^{-1} \, q),
\end{equation}  
so that {\rm (\ref{E:3.30})} and {\rm (\ref{E:3.11aa})} are valid, 
and  $F$ is  skew-symmetric with entries in ${\rm BMO}(\R^n)$. 

{\rm(iii)} If $n=2$ then the bilinear form of $\mathcal{L}$  is bounded  
if and only if   
${\rm div} \, \vec b = q = 0$, and $\vec b = (\partial_2 g, -\partial_1 g)$, 
where $g = \Delta^{-1} {\rm rot} \, \vec b \in {\rm BMO}(\R^2)$. 
\end{theorem} 

\begin{proof} Let us first prove the sufficiency part of statement (i). 
To see that 
$q = {\rm div} \, \vec h$ is form bounded 
on $L^{1, \, 2} (\R^n)$, we 
use integration by parts, Schwarz's inequality, and (\ref{E:3.11aa}):
\begin{align}
| \langle q u, \, v\rangle | & 
= | \langle \vec h, \nabla \, 
( \bar u \, v\rangle) | =  | \langle \vec h, \,  v  \, \nabla \bar u +  
 \bar u \, \nabla v \rangle) |   \notag \\ 
& \le C \, ( ||u \,  \vec h||_{L^2(\R^n)} \, ||\nabla \bar v||_{L^2(\R^n)}
+ ||\bar v \,  \vec h||_{L^2(\R^n)} \, ||\nabla u||_{L^2(\R^n)} ) \notag 
\\ & \le C \, ||\nabla u||_{L^2(\R^n)} \, ||\nabla v||_{L^2(\R^n)},
\notag
\end{align}
for every $u$, $v \in C^\infty_0(\R^n)$.

Suppose now that  $\vec b$ is represented in the form {\rm (\ref{E:3.30})}, 
i.e., 
$\vec b = \vec c + \vec d,$ where 
 $\vec c$ satisfies   inequality 
{\rm (\ref{E:3.11aa})}, and $\vec d = {\rm Div} \, F$.   
Then  by Schwarz's inequality 
and {\rm (\ref{E:3.11aa})}, it follows: 
$$
| \langle \vec c \cdot \nabla u, \, v \rangle|  
\le \, ||v \, \vec c||_{L^2(\R^n)} \, ||\nabla u||_{L^2(\R^n)} 
\le C \, \, ||\nabla u||_{L^2(\R^n)} \, ||\nabla v||_{L^2(\R^n)},  
$$
which proves that $\vec c \cdot \nabla$ is form bounded. 

We next treat $\vec d \cdot \nabla$. 
By Proposition~\ref{Proposition 3.1} in the case $q=0$,   
this is equivalent to both {\rm (\ref{E:3.2})} and {\rm (\ref{E:3.3})} 
with $\vec d$ in place of $\vec b$. 
Since ${\rm div} \, \vec d = 0$, the first 
condition becomes vacuous, and so it suffices to verify  
(\ref{E:2.14}). The latter inequality holds by Lemma~\ref{Lemma 3.8} 
since $\vec d = {\rm Div} \, F$ where 
$F \in {\rm BMO}(\R^n)^{n \times n}$, and $F$ is skew-symmetric. 
Combining 
the preceding estimates we conclude that   $\mathcal{L} = 
\vec c \cdot \nabla + \vec d \cdot \nabla +q$   
is form bounded.

Conversely, if the bilinear form of $\mathcal{L}$ is bounded, then by 
Lemma~\ref{Lemma 3.3} and Lemma~\ref{Lemma 3.5},  
decomposition (\ref{E:3.30}) holds, 
where $\vec c$, $F$, and $\vec h$ are given by (\ref{E:3.delta}). 
Furthermore,  it follows that 
 (\ref{E:3.11aa}) holds for $n \ge 3$, and 
 (iii) is valid if $n=2$. The proof of Theorem~\ref{Theorem 3.3} 
is complete. \end{proof}

\begin{corollary}\label{Corollary 2.2}  
Let $n \ge 3$. 
Let $\mathcal{L} =  \vec b \, \cdot \nabla + q$, where 
$\vec b \in D'(\R^n)^n$, and $q \in D'(\R^n)$. 
 Then  
the bilinear form of $\mathcal{L}$ is bounded on 
$L^{1, \, 2} (\R^n) \times L^{1, \, 2} (\R^n)$
 if and only if 
 decomposition {\rm (\ref{E:3.30})} is valid where 
$\vec c$, $F$ and $\vec h$ are given by {\rm (\ref{E:3.delta})}, 
$F = \Delta^{-1} {\rm curl} \, \vec b 
\in {\rm BMO}(\R^n)^{n \times n}$,  and the locally 
finite measure 
$$d \mu(x) = \left ( |\nabla (\Delta^{-1}  q)|^2 +   
|\nabla ( \Delta^{-1} {\rm div} \, \vec b) |^2 \right) \, dx $$ 
 is subject to any one of the following equivalent conditions:

 {\rm (i)} For every compact set $e\subset \R^n$, 
$$\mu (e) \le C \,  {\rm cap} \, (e).$$ 
where the capacity  ${\rm cap} \, (\cdot)$ 
 is defined by {\rm (\ref{E:cap})}. 

{\rm (ii)} For any cube $P$ in $\R^n$, 
$$\int_P \left[ (-\Delta)^{-\frac 1 2} 
(\chi_P \, d \mu) \right]^2 \, dx \le C \, \mu(P),
$$
where $C$ does not depend on $P$.

{\rm (iii)} For a.e. $x \in \R^n$,
$$ (-\Delta)^{-\frac 1 2} 
\left [ (-\Delta)^{-\frac 1 2} \mu \right]^2 (x) \le C \, 
(-\Delta)^{-\frac 1 2} \mu (x) < + \infty.
$$

{\rm (iv)} For any dyadic cube $P$ in $\R^n$,
$$\sum_{Q \subseteq P} \, \frac {\mu(Q)^2}  
{|Q|^{1- \frac 2 n}} \le C \, \mu(P),$$
where the sum is taken over all dyadic cubes $Q$ 
contained in $P$, and $C$ does not depend on $P$.
\end{corollary}

Corollary~\ref{Corollary 2.2}  follows by coupling Theorem~\ref{Theorem 3.3} 
with Theorem~\ref{Trace Theorem}.

\begin{corollary} \label{Corollary 2.3} 
Let $n \ge 3$. 
Let $\mathcal{L} = \vec b \, \cdot \nabla + q$, 
where $\vec b \in D'(\R^n)^n$, and $q \in D'(\R^n)$.

{\rm (i)} If 
the sesquilinear form of $\mathcal{L}$ is bounded on 
$L^{1, \, 2} (\R^n) \times L^{1, \, 2} (\R^n)$, then 
{\rm (\ref{E:3.30})} and {\rm (\ref{E:3.delta})} hold 
and, for every ball $B$ in $\R^n$,  
\begin{equation}\label{E:3.12}
\int_{B} [ \, |\nabla (\Delta^{-1} {\rm div} \, \vec b) |^2 + 
|\nabla \Delta^{-1} \, q|^2 ] \, dx \le C \, |B|^{1 - \frac 2 n}. 
\end{equation}

{\rm (ii)}  Conversely, if  {\rm (\ref{E:3.30})} 
and {\rm (\ref{E:3.delta})} 
are valid where 
\begin{equation}\label{E:3.13}
\int_{B} [ \, |\nabla ( \Delta^{-1} {\rm div} \, \vec b)|^{2} + 
|\nabla \Delta^{-1} \, q|^{2} ]^{1 + \epsilon}
 \, dx \le C \, |B|^{1- \frac{2 (1 + \epsilon)} n}, 
\end{equation} 
for some $\epsilon >0$, and $\Delta^{-1} {\rm curl} \, \vec b 
\in {\rm BMO}(\R^n)^{n \times n}$,  
then  the sesquilinear form of $\mathcal{L}$ is  bounded on 
$L^{1, \, 2} (\R^n) \times L^{1, \, 2} (\R^n)$. 
\end{corollary}

\begin{remark} Inequality (\ref{E:3.12}), together 
with  Poincar\'e's inequality, yields that 
$\Delta^{-1} {\rm div} \, \vec b\in {\rm BMO}(\R^n)$ 
and  $\Delta^{-1} q \in {\rm BMO}(\R^n)$.  
\end{remark}

Statement (i) of Corollary~\ref{Corollary 2.3} 
is immediate from Theorem~\ref{Theorem 3.3} and
(\ref{E:cube}), whereas (ii) follows by combining 
Corollary~\ref{Corollary 2.2} with the 
Fefferman--Phong condition (\ref{E:1.5b}) applied to 
$d \mu(x) = \left( |\nabla (\Delta^{-1}  {\rm div} \, \vec b) |^2 + 
|\nabla (\Delta^{-1} q) |^2 \right) 
\, dx $. Sharper sufficient conditions 
are deduced from  Corollary~\ref{Corollary 2.2}  in the same way 
by making use of  
the conditions  due to Chang, Wilson, and Wolff \cite{ChWW} applied to 
$d \mu$.

The following statement 
is a consequence of  Theorem~\ref{Theorem 3.3}   
and Lemma~\ref{Lemma 3.3}. 
\begin{corollary}\label{corollary 2.4} Let $\vec b \in D'(\R^n)^n$ and 
$q \in D'(\R^n)$, $n \ge 2$.  Then 
the operator 
$\vec b \cdot \nabla +q$ is form bounded on $L^{1, \, 2}(\R^n)\times L^{1, \, 2}(\R^n)$ if and 
only if, for every cube $Q$ in $\R^n$, 
\begin{equation}\label{E:1.loc}
||\vec b||_{L^{1, \, 2}(Q)} \le {\rm const} \, |Q|^{\frac  1 2}
\end{equation}
and both 
\begin{equation}\label{E:1.mult}
{\rm div} \, \vec b : L^{1, \, 2}(\R^n) \to L^{-1, \, 2}(\R^n),
\quad q : L^{1, \, 2}(\R^n) \to L^{-1, \, 2}(\R^n),
\end{equation}
are bounded multiplication operators.
\end{corollary}

We conclude this section with a form boundedness 
criterion for the magnetic Schr\"odinger operator
 $$\mathcal{M} = (i \, \nabla +\vec a )^2 + q$$
with real-valued magnetic vector 
potential $\vec a$. As a direct consequence of Theorem~\ref{Theorem 3.3}, 
we establish the following form boundedness criterion for 
$\mathcal{M}$.

\begin{theorem}\label{Theorem 3.4}
Let  $\vec a \in L^2_{{\rm loc}}(\R^n)$ 
and 
$q\in D'(\R^n)$, $n \ge 2$. Then the operator 
$\mathcal{M}=(i \, \nabla +\vec a )^2 + q$  is form 
bounded on $L^{1, \, 2} (\R^n)\times L^{1, \, 2} (\R^n)$ 
if and only if both $q + |\vec a|^2$ and $\vec a \cdot \nabla$ 
are form bounded. More precisely, in order that 
\begin{equation}
| \langle \mathcal{M} \, u, \, v \rangle | \le C \, 
||u||_{L^{1, \, 2} (\R^n)} \, ||v||_{L^{1, \, 2} (\R^n)}, 
\quad   u, \, v \in C^\infty_0(\R^n), 
\end{equation}
it is necessary and sufficient that  
\begin{equation}\label{E:3.14}
\vec a = \vec c + {\rm Div} \, F, \qquad q + |\vec a|^2 = 
{\rm div} \, \vec h, 
 \end{equation}
where $F$ is a skew-symmetric matrix field whose entries belong to 
${\rm BMO}(\R^n)$, and $|\vec c|^2 + |\vec h|^2   \in \mathfrak{M}^{1, \, 2}_+ (\R^n)$. 

Moreover, one can define $\vec c $, $F$, and $\vec h$ in  representation 
{\rm (\ref{E:3.14})}  
constructively as 
 $\vec c = \nabla ( \Delta^{-1} {\rm div} \, \vec a)$, $F = 
 \Delta^{-1} \, {\rm curl} \, \vec a$, 
and $\vec h = \nabla \Delta^{-1} \, (q + |\vec a|^2)$. 

In the case $n=2$, $\mathcal{M}$ is form bounded on $L^{1, \, 2} (\R^2)\times L^{1, \, 2} (\R^2)$
if and only if 
${\rm div} \, \vec a =0$, and $q + |\vec a|^2 =0$, where $\vec a = 
( \partial_2 g, \, -\partial_1 g)$, and $g\in {\rm BMO} (\R^2)$.  
\end{theorem}

\begin{remark} This characterization simplifies  
under the Coulomb  gauge hypothesis 
${\rm div} \, \vec a =0$ (see \cite{RS}, Sec. X. 4). Then, for the form 
boundedness of $(i \, \nabla + \vec a)^2 + q$ on $L^{1, \, 2} (\R^n)\times L^{1, \, 2} (\R^n)$,  
$n \ge 3$, 
it is necessary 
and sufficient that $q + |\vec a|^2$ be form bounded, and 
$\vec a = {\rm Div} \, F$, 
where $F = \Delta^{-1} \, {\rm curl} \, \vec a \in 
{\rm BMO} (\R^n)^{n \times n}$.\end{remark}

\begin{remark} The above characterization 
of the form boundedness of $\mathcal{M}$ 
 holds  if one replaces the assumption 
$\vec a \in L^2_{{\rm loc}}(\R^n)$ by 
$q+ |\vec a|^2 \in L^1_{{\rm loc}}(\R^n)$. \end{remark}

\section{An estimate for $\langle \, | \vec b\cdot \nabla u|,  \, u \, 
\rangle$}  \label{Section 4} 

In this section we prove the following statement 
for the nonlinear quadratic form 
$\langle \, | \vec b\cdot \nabla u|,  \, u \, 
\rangle$ 
which holds for every 
open set $\Omega\subset \R^n$, in particular $\Omega = \R^n$.

\begin{proposition}\label{Proposition 4.1} 
 Let $\vec b \in L^1_{{\rm loc}} (\Omega)^n$. Then the  
best constant in the inequality 
\begin{equation}\label{E:4.1}
\left | \int_\Omega | \vec b \cdot \nabla u| \, \bar u  \, dx \right | \le C \, 
||\nabla u||^2_{L^2 (\Omega)}, \qquad   u \in C^\infty_0(\Omega),
\end{equation} 
satisfies the estimates 
\begin{equation}\label{E:4.2}
C \le c \le 2 \sqrt{n} \, C,
\end{equation} 
where $c^2$ is the best constant in the  inequality 
\begin{equation}\label{E:4.3}
\int_\Omega |\vec b|^2  \, |u|^2 \, dx \le c^2 \, 
||\nabla u||^2_{L^2 (\Omega)}, \qquad   u \in C^\infty_0(\Omega). 
\end{equation} 
\end{proposition} 

\begin{remark} The constant $c$ in the previous inequality coincides 
with the norm of the multiplier operator 
$\vec b : \, L^{1, \, 2} (\Omega) \to L^{2} (\Omega)^n$  
where $L^{1, \, 2} (\Omega)$ is a (homogeneous) Sobolev space defined as 
the completion of $C^\infty_0(\Omega)$ in the Dirichlet norm 
$||\nabla u||_{L^2 (\Omega)}$.
\end{remark}

\begin{proof}The lower bound in (\ref{E:4.2}) is obvious. 
Let us prove the upper bound. For real-valued $u$, inequality (\ref{E:4.1}), combined with 
the well-known estimate $|| \nabla |u| \,  ||_{L^2(\Omega)} \le 
|| \nabla u ||_{L^2(\Omega)}$  (see \cite{LL}, Sec. 7.8), yields 
\begin{equation}\label{E:4.1a}
 \int_\Omega | \vec b \cdot \nabla u^2 |  \, dx \le C \, 
||\nabla u||^2_{L^2 (\Omega)}, \qquad   u \in C^\infty_0(\Omega). 
\end{equation} 
 Consequently, 
for complex-valued $u$, 
\begin{equation}\label{E:4.1b}
 \int_\Omega | \vec b \cdot \nabla u^2 |  \, dx \le 2 C \, 
||\nabla u||^2_{L^2 (\Omega)}, \qquad   u \in C^\infty_0(\Omega). 
\end{equation} 
For every $\epsilon>0$, by (\ref{E:4.1b}),
$$\int_\Omega |( \vec b \cdot \nabla) (\epsilon \, u \pm  
\epsilon^{-1} \, v)^2|   \, dx \le 2 C \, ||\epsilon \, \nabla u \pm  
\epsilon^{-1} \, v ||^2_{L^2 (\Omega)}.
$$
Hence,
\begin{align}
4 \, \int_\Omega  |( \vec b \cdot \nabla) (u \, v)| \, dx & = 
\int_\Omega  |( \vec b \cdot \nabla) \left ( (\epsilon \, u +
\epsilon^{-1} \, v)^2 -  (\epsilon \, u - 
\epsilon^{-1} \, v)^2 \right)| \, dx \notag \\ & \le 4C \, 
\left (\epsilon^2 \, ||\nabla u||^2_{L^2 (\Omega)} + 
\epsilon^{-2} \, ||\nabla v||^2_{L^2 (\Omega)} \right).
\notag
\end{align}
Minimizing over $\epsilon$, we get:
$$
 \int_\Omega  |( \vec b \cdot \nabla) (u \, v)| \, dx \le 2 C \, 
||\nabla u||_{L^2 (\Omega)} \, ||\nabla v||_{L^2 (\Omega)}.
$$

We now set 
$$
u(x) = e^{i \, \langle \omega, \, \xi\rangle} \, h(x), 
\qquad h \in C^\infty_0(\Omega), 
$$
where $\omega \in S^{n-1}$ and $\xi \in \R^n$. We estimate:
$$
\int_\Omega |\langle \omega, \, \vec b \rangle| \, |h \, \, v| \, dx \le 
2 C \, ||h||_{L^2 (\Omega)} \, ||\nabla v||_{L^2 (\Omega)} + O (|\xi|^{-1}).
$$
Letting $|\xi|\to +\infty$ gives 
$$
\int_\Omega |\langle \omega, \, \vec b \rangle|^2 \, |v|^2 \, dx 
\le 4 C^2 \, 
||\nabla v||_{L^2 (\Omega)}.
$$
Integrating the preceding inequality over $S^{n-1}$ and using the identity 
$$\int_{S^{n-1}} |\langle \omega, \, \vec b \rangle|^2 \, d s_\omega = 
\frac 1 n |S^{n-1}| \, |\vec b|^2,
$$
we arrive at:
$$
\int_\Omega |\vec b|^2 \, |v|^2 \, dx \le 4 n C^2 \, 
||\nabla v||_{L^2 (\Omega)}.
$$
The proof of Proposition~\ref{Proposition 4.1} is complete. 
\end{proof}

\begin{corollary}\label{Corollary 4.2} 
 The best constant in the inequality {\rm (\ref{E:4.1})} satisfies 
the estimates: 
\begin{equation}\label{E:4.5}
\tfrac 1 2 \, C \le \sup_{e\subset \Omega} \frac { ||\vec b||_{L^2(e)} } 
{{\rm cap} \, (e, \, \Omega)^{\frac 1 2}} 
\le 2 \sqrt{n} \, C,
\end{equation}
where the supremum is taken over all compact sets  
$e \subset \Omega$ of 
positive capacity defined by 
$$ {\rm cap} \, (e, \, \Omega) = 
\inf \, \left\{ \, ||\nabla u||^2_{L^2(\Omega)} \, 
: \quad u \in C^\infty_0(\Omega),  \quad u \ge 1 \, \, {\rm on} \,
 \, e \right\}.
$$
\end{corollary}

Corollary~\ref{Corollary 4.2} follows from 
Proposition~\ref{Proposition 4.1} 
and \cite{M}, Sec. 2.5.

\section{Form boundedness on the Sobolev space $W^{1, \, 2} (\R^n)$}
\label{Section 5} 

In this section, we obtain the form boundedness 
criterion for the general second order differential operator $\mathcal{L}$ on 
the  Sobolev space $W^{1, \, 2} (\R^n)$, $n \ge 2$. As was noticed above, without 
loss of generality we may assume that $\mathcal{L}$ is in the divergence form:  
$\mathcal{L} = {\rm div} \, (A \, \nabla) +   \vec b \cdot \nabla +q$, where  
 $A \in D'(\R^n)^{n \times n}$, 
$\vec b \in D'(\R^n)^n$ and $q \in D'(\R^n)$, $n \ge 2$. 
Then clearly the sesquilinear inequality 
\begin{equation}\label{E:5.1}
\left \vert \langle \mathcal{L} \, u, \, v  \rangle\right \vert 
\le C \,  ||u||_{W^{1, \, 2} (\R^n)} ||v||_{W^{1, \, 2} 
(\R^n)}
\end{equation} 
holds for all $u, \, v \in C^\infty_0(\R^n)$ if and only if  the operator $\mathcal{L}$ 
(or, more precisely, its unique 
extension from $C^\infty_0$ to $W^{1, \, 2}$), 
\begin{equation}\label{E:5.2}
\mathcal{L} \, : \, W^{1, \, 2} (\R^n) 
\to W^{-1, \, 2} (\R^n) 
\end{equation}
is bounded.  

We notice that Proposition~\ref{Proposition 2.100} 
holds for $W^{1, \, 2}$   in place of $L^{1, \, 2}$, with obvious modifications in the proof. In particular, the condition 
$A^s = \tfrac 1 2 (A+A^t)  \in L^\infty(\R^n)^{n \times n}$ is necessary 
for the form boundedness of $\mathcal{L}$ on $W^{1, \, 2}$, whereas  $A^c = \tfrac 1 2 (A-A^t)$ 
can be included in $\vec b$ by letting $\vec b_1 = \vec b - {\rm Div} \, A^c$, 
exactly as in the case of the homogeneous Sobolev space. In other words, 
it suffices to consider the form boundedness problem on $W^{1, \, 2}(\R^n)$ 
 for ${\mathcal L} = \vec b \cdot \nabla +q$.

Recall that ${\rm BMO}^\#(\R^n)$ stands for the space of 
$f \in L^1_{{\rm loc}}(\R^n)$ 
such that 
$$
\sup_{x_0 \in \R^n, \, 0<\delta\le 1} \, \frac 1 {|B_\delta(x_0)|} 
\int_{B_\delta(x_0)} |f(x) - m_{B_\delta(x_0)}(f)| \, dx < +\infty.
$$

\begin{theorem}\label{Theorem 5.1}
Let $\vec b \in D'(\R^n)^n$, $q \in D'(\R^n)$,  and  let 
$\mathcal{L} =  \vec b \cdot \nabla +q$, $n \ge 2$. 
Then {\rm (\ref{E:5.1})} holds if and only if $\vec b$ and $q$ 
can be represented 
respectively in 
the form:
\begin{equation}\label{E:5.3}
\vec b = \vec c + {\rm Div} \, F, 
\quad q = {\rm div} \, \vec h + \gamma,
\end{equation} 
where $F$ is a skew-symmetric matrix field such that 
$F \in  {\rm BMO}^\#(\R^n)^{n \times n}$, and   $(|\vec c|^2 + |\vec h|^2 + |\gamma|) \, dx$ 
is an admissible measure for $W^{1, \, 2} (\R^n)$, i.e., 
\begin{equation}\label{E:5.4}
\int_{\R^n} (|\vec c|^2 + |\vec h|^2 + |\gamma|) \, 
|u|^2 \, dx  \le c \, 
||u||^2_{W^{1, \, 2} (\R^n)}, \qquad   u \in C^\infty_0(\R^n).
\end{equation}
Moreover, in the decomposition {\rm (\ref{E:5.3})} and 
condition {\rm (\ref{E:5.4})}
one can set 
\begin{align}\label{E:5.5a}
& \vec c = - \nabla [(1- \Delta)^{-1} {\rm div} \, \vec b] + 
(1 - \Delta)^{-1} \, 
\vec b, 
\qquad  F =   - (1-\Delta)^{-1} {\rm curl} \, \vec b, \\
& \vec h = - \nabla (1- \Delta)^{-1} q, \qquad 
\gamma = (1- \Delta)^{-1} q. 
\label{E:5.5b}
\end{align} 
Furthermore, 
{\rm (\ref{E:5.4})} holds with 
$|(1- \Delta)^{-1} {\rm div} \, \vec b|^2
+ |(1- \Delta)^{-1} \, \vec b|^2$ in place of $|\vec c|^2$. 
\end{theorem}

\begin{proof} Suppose that $\vec b$ is given by (\ref{E:5.3}) 
where $F$ is a skew-symmetric matrix field such that 
$F \in  {\rm BMO}^\#(\R^n)^{n \times n}$, and   
$\vec c$, $\vec h$, and $\gamma$ satisfy (\ref{E:5.4}). The boundedness 
of the bilinear  form associated with $q$ and $\vec c\cdot\nabla$  
follows easily using 
 integration by parts and Schwarz's inequality:
\begin{align}
 & \left \vert \langle \vec c \cdot \nabla u  + 
q \, u, \, v  \rangle \right \vert 
 \le \left \vert \langle \vec c \cdot \nabla u,  \, v  \rangle\right \vert
+  \left \vert \langle \vec h, \, \bar u \, \nabla v 
+ v \, \nabla \bar u \rangle\right \vert + 
\left \vert \langle \vec \gamma, \, \bar u \,  v \rangle\right \vert 
\notag \\ & \le \, \left \Vert |\vec c| \, |v| \right\Vert_{L^2(\R^n)} \, 
||\nabla u||_{L^2(\R^n)} + 
\left \Vert |\vec h| \, |v| \right\Vert_{L^2(\R^n)} \, 
||\nabla u||_{L^2(\R^n)} \notag \\& + 
\left \Vert |\vec h| \, |u| \right\Vert_{L^2(\R^n)} \, 
||\nabla v||_{L^2(\R^n)} + 
\left \Vert |\vec \gamma|^{\frac 1 2} \, |u| \right\Vert_{L^2(\R^n)} \, 
\, \left \Vert |\vec \gamma|^{\frac 1 2} \, |v| \right\Vert_{L^2(\R^n)} 
 \notag \\ 
& \le C \,  ||u||_{W^{1, \, 2} (\R^n)} \, 
||v||_{W^{1, \, 2} (\R^n)}, \qquad   \, 
u, \, v \in C^\infty_0(\R^n). \notag
\end{align}

We next prove the boundedness of the bilinear form associated with 
the divergence free part of $\vec b$ given by 
$\vec d = {\rm Div} \, F$. This may be  viewed as an inhomogeneous 
version of the  ${\rm div}$-${\rm curl}$ lemma \cite{CLMS}. 
The proof is based on a localization 
principle, combined with an appropriate extension of ${\rm BMO}(B)$ 
functions originally defined on a ball $B\subset \R^n$. 
 
\begin{lemma}\label{Lemma 5.2} Suppose  $\vec d = {\rm Div} \, F$ 
in $D'(\R^n)^n$, where $F$ is a skew-symmetric matrix 
function such that $F \in {\rm BMO}^\# (\R^n)^{n\times n}$.
Then the inequality
\begin{equation}\label{E:5.7}
\left \vert \langle \vec d \cdot \nabla u, \, v  \rangle\right \vert 
 \le C \,  ||u||_{W^{1, \, 2} (\R^n)} \, ||v||_{W^{1, \, 2} (\R^n)},   
\qquad   \, 
u, \, v \in C^\infty_0(\R^n),   
\end{equation}
holds where $C$ does not depend on $u$ and $v$. 
\end{lemma}

\begin{proof} We first prove a localized version of 
{\rm (\ref{E:5.7})},  
  \begin{equation}\label{E:5.9}
\left \vert \langle \vec d \cdot \nabla u, \, v  \rangle\right \vert 
 \le C \,  ||\nabla u||_{L^2 (B_1(x_0))} \, 
||\nabla v||_{L^2 (B_1(x_0))}, 
\end{equation}
where the constant $C$ does not depend on $u, \, v 
\in C^\infty_0(B_1(x_0))$, and $x_0 \in \R^n$.

For a domain $\Omega \subset \R^n$, denote by ${\rm BMO}(\Omega)$ 
 the space of functions $f\in L^1_{{\rm loc}} (\Omega)$ 
 such that 
$$
\sup_{B\subset \Omega} \frac 1 {|B|} \int_B |f - m_B(f)| \, dx < +\infty,
$$
where the supremum is taken over all balls $B$ in $\Omega$. 

Since 
$F \in {\rm BMO}^\#(\R^n)^{n\times n}$, 
it follows that
  $F \in {\rm BMO}(B_1(x_0))$ for every $x_0 \in \R^n$, and 
\begin{equation}\label{E:5.10}
\sup_{x_0\in \R^n} \, ||F||_{{\rm BMO}(B_1(x_0))^{n\times n}} 
< +\infty.
\end{equation}

By replacing $u$ and $v\in C^\infty_0(B_1(x_0))$ in {\rm (\ref{E:5.9})} 
with  $u(x-x_0)$ and $v(x-x_0)$ respectively, one 
can assume without loss of 
generality that $x_0=0$, and $F \in{\rm BMO}(B_1(0))$.  
Denote by 
$\widetilde F$ an extension of $F$ from $B_1(0)$
to $\R^n$ such that 
\begin{equation}\label{E:5.11}
||\widetilde F||_{{\rm BMO}(\R^n)^{n\times n}} \le c \, 
||F||_{{\rm BMO}(B_1(0))^{n\times n}},
 \end{equation}
where $c$ depends only on $n$. To construct such an extension one can 
use a reflection 
in the boundary. (See, e.g., \cite{J} where this is done for very general
 domains $\Omega\subset \R^n$.) 

Note that both $F$ and 
$\widetilde F$ are 
skew-symmetric. Hence by (\ref{E:5.11}) and 
the  version of the ${\rm div}$-${\rm curl}$ lemma 
used above (see the proof of 
Lemma~\ref{Lemma 3.8}), 
\begin{align}
2 \, \left \vert \langle \vec d \cdot \nabla u, \, v  
\rangle\right \vert 
 & =  \, \left \vert \int_{B_1(0)} {\rm trace} \, 
F \cdot  
\{ \partial_i \bar u \, \partial_j v - \partial_j \bar u \, 
\partial_i \bar v\} \, 
dx \right \vert \notag \\
& =  \, \left \vert \int_{B_1(0)} {\rm trace} \, 
\widetilde F \cdot
\{\partial_i \bar u \, \partial_j v - \partial_j \bar u \, 
\partial_i v\}  \, 
dx \right \vert \notag \\ 
& \le C \, || F||_{{\rm BMO} (\R^n)^{n \times n}}  \,  
 ||\nabla u||_{L^2(B_1(0))} \, ||\nabla v||_{L^2(B_1(0))}, 
\notag
\end{align}
where $C$ depends only on $n$. Taking into account 
{\rm (\ref{E:5.10})}, we conclude that 
{\rm (\ref{E:5.9})} holds 
for every $u, \, v \in C^\infty_0(B_1(x_0))$ with a constant which 
does not depend on $u, \, v$, and $x_0$.

To prove {\rm (\ref{E:5.7})}, suppose 
$u, \, v \in C^\infty_0(B_R(x_0))$, 
$R>1$.  Pick a sequence of functions 
$\{\zeta_i\}_{i=1}^\infty$  so that 
 \begin{align}
& \sum_i \, \zeta_i(x)^2 =1,   
\quad \sum_i \, |\nabla \zeta_i(x)|^2 \le c(n) \quad {\rm on} \, \, 
B_R(x_0), \label{E:5.12} \\ & \sum_i \, \zeta_i^2 \in C^\infty(\R^n), 
\quad 
\zeta_i \in C^\infty_0(B_1(x_i)), \quad 
i =1, 2, \ldots. \label{E:5.13}
\end{align}
Here $x_i$ is a cubic lattice of equidistant 
points in $\R^n$ with grid distance equal to $\frac 1 {2 \sqrt{n}}$. 
(See, e.g.,  \cite{MV4}, the proof of Lemma 3.1). 

Now integration by parts gives
\begin{align}
 \langle \vec d \cdot \nabla u, \, v  \rangle  & = 
\sum_i \, \langle \vec d \cdot \nabla u, \, \zeta_i^2 \, v 
 \rangle\notag 
\\ & = \sum_i \, \langle \vec d \cdot \nabla (\zeta_i \, u), \, 
\zeta_i\, v  \rangle
- \tfrac 1 2 \, \sum_i \, \langle \vec d \cdot 
 \nabla (\zeta_i^2), \bar u \, v  \rangle\notag \\
& =  \sum_i \, \langle \vec d \cdot \nabla (\zeta_i \, u), \, 
\zeta_i\, v  \rangle. 
\notag
\end{align}
In the last line we have used 
$\sum_i \, \nabla (\zeta_i^2) =0$ on $B_R(x_0)$ which 
follows from (\ref{E:5.12}).

Suppose now that (\ref{E:5.9}) holds. 
Then from the preceding equation  we deduce: 
\begin{align}
\left \vert \langle \vec d \cdot \nabla u, \, v  \rangle \right 
\vert  & \le C \, \sum_i \, 
\left \vert \langle \vec d \cdot \nabla (\zeta_i u), \, 
\zeta_i \, v  \rangle \right \vert \notag\\
& \le C \, \sum_i \, ||\nabla \, (\zeta_i \, u)||_{L^2 (\R^n)} 
||\nabla \, (\zeta_i \, v)||_{L^2 (\R^n)} \notag \\ & \le 
C \, \sum_i \, ||\nabla \, (\zeta_i \, u)||^2_{L^2 (\R^n)} + 
C \, \sum_i \, ||\nabla \, (\zeta_i \, v)||^2_{L^2 (\R^n)}. 
\notag
\end{align}
We estimate the first term on the 
right-hand side using  (\ref{E:5.12}):
\begin{align}
  \sum_i \, ||\nabla \, (\zeta_i \, u)||^2_{L^2 (\R^n)} 
& \le \,  C \, \sum_i \,  ||\zeta_i \, \nabla u||^2_{L^2 (\R^n)}
+ C \sum_i \, ||(\nabla \zeta_i) \, u||^2_{L^2 (\R^n)} \notag \\ 
& \le C \, ||u||^2_{W^{1, \, 2} (\R^n)}. 
\notag
\end{align}
where $C$ does not depend on $u$. A similar estimate holds 
for the second term which involves $v$.  Note that without loss of generality 
we may assume that  $\max \left( || u||_{W^{1, \, 2} (\R^n)}, \, 
|| v||_{W^{1, \, 2} (\R^n)}    \right) \le 1$. 
Hence, 
$
\left \vert \langle \vec d \cdot \nabla u, \, v  \rangle \right 
\vert \le C,$ 
which yields (\ref{E:5.7}). This concludes the proof of 
Lemma~\ref{Lemma 5.2}. 
\end{proof} 

It follows from the preceding estimates for $q$ and $\vec c$, 
and Lemma~\ref{Lemma 5.2} that 
$\mathcal{L}=\vec b \cdot \nabla +q $ is form bounded provided  (\ref{E:5.3}) holds 
with $F \in {\rm BMO}^\#(\R^n)^{n \times n}$ and $\vec c$, $\vec h$, $\gamma$ 
satisfying (\ref{E:5.4}).

It remains to prove the converse for $F $, $\vec c$, $\vec h$, and 
$\gamma$ 
defined by (\ref{E:5.5a})--(\ref{E:5.5b}). It is easy to see that 
Proposition~\ref{Proposition 3.1} 
holds verbatim with 
$W^{1, \, 2}$  
in place of $L^{1, \, 2}$; i.e., (\ref{E:5.1}) holds if and only 
both of the following inequalities are valid:
\begin{align}\label{E:5.1a}
\left \vert
 \langle (q - \tfrac 1 2 \, {\rm div} \, \vec b)
\, u, \, v \rangle \right \vert & \le C \,
||u||_{W^{1, \, 2}(\R^n)} \, ||v||_{W^{1, \, 2}(\R^n)}, \\ 
\left \vert
 \langle  \vec b,  \,
 \bar u \, \nabla v  - v \nabla \bar u \rangle \right \vert
& \le C \,
||u||_{W^{1, \, 2}(\R^n)} \, ||v||_{W^{1, \, 2}(\R^n)},
\label{E:5.1b}
\end{align} 
for all $u, \, v \in C^\infty_0(\R^n)$.

 An analogue of 
Proposition~\ref{Proposition 3.2} states that, if  
 {\rm (\ref{E:5.1b})} holds, 
then the following estimates are valid:
\begin{align}\label{E:5.2a}
|| {\rm div} \, \vec b||_{W^{-1, \, 2}(Q)} & \le C \,
|Q|^{\frac 1 2 - \frac 1 n} \quad {\rm if} \, \, n \ge 3, 
\\  || {\rm div} \, \vec b||_{W^{-1, \, 2}(Q)} & \le C \,
\left ( \log \, \tfrac 2 {|Q|}\right )^{-\frac 1 2}  \quad {\rm if} \, \, n =2,
\label{E:5.2b}\\ ||\vec b||_{W^{-1, \, 2}(Q)} & \le C \, |Q|^{\frac 1 2}
\quad {\rm if} \, \, n \ge 2,
\label{E:5.2c}
\end{align}
for every cube $Q$ in $\R^n$ such that $\ell(Q) \le 1$. The only 
change that is needed in the proof is that, for the capacity 
${\rm Cap} \, (\cdot)$ associated with $W^{1, \, 2}(\R^n)$,  which is defined by 
(\ref{E:Cap}), 
we have 
${\rm Cap} \, (Q)  \simeq  \left ( \log \, \tfrac 2 {|Q|}\right )^{-1}$
for $n \ge 2$ and $\ell(Q) \le 1$ by  (\ref{E:cube}). 
(Note that in two dimensions, contrary to the case of  
$L^{1, \, 2}(\R^2)$, $\vec b$  is no longer required to 
be  divergence free.)

It now follows from {\rm (\ref{E:5.2c})}, as in the proofs of Lemma~\ref{Lemma 3.3} and  ~\ref{Lemma 3.5}, that decomposition 
{\rm (\ref{E:5.3})} holds where $\vec c$, $F$, $\vec h$, 
and $\gamma$ are given by {\rm (\ref{E:5.5a})} and {\rm (\ref{E:5.5b})} 
respectively, and $F \in {\rm BMO}^\#(\R^n)^{n \times n}$. Furthermore, using a 
direct analogue 
of  Lemma~\ref{Lemma 3.6}  for $W^{1, \, 2}(\R^n)$, we deduce from 
(\ref{E:5.1b}) that ${\rm div} \, \vec b$ is form bounded on  $W^{1, \, 2}(\R^n)$, i.e., 
$$| \langle ( {\rm div} \, \vec b) \, u, \, v\rangle | \le C \, ||u||_{W^{1, \, 2} (\R^n)} \, 
||v||_{W^{1, \, 2} (\R^n)}, 
$$
for all $u, v \in C^\infty_0(\R^n)$. Hence, by  (\ref{E:5.1a}),  
$$| \langle q \, u, \, v\rangle | \le C \, ||u||_{W^{1, \, 2} (\R^n)} \, 
||v||_{W^{1, \, 2} (\R^n)}, 
$$
for all $u, v \in C^\infty_0(\R^n)$.

The preceding inequality,  by  Theorem 4.2 in \cite{MV1},  yields
\begin{equation}\label{E:5.14} 
\int_{\R^n}  ( | \nabla (1-\Delta)^{-1} q |^2 + |(1-\Delta)^{-1} q| ) \, |u|^2 \, dx \le C \, 
||u||^2_{W^{1, \, 2}(\R^n)}, 
\end{equation}
for all $u \in C^\infty_0(\R^n)$. Note that,  according to \cite{MV1} (Sec. 4, Remark 3),  
 it is possible to put  $|(1-\Delta)^{-1} q |^2$ in place of $|(1-\Delta)^{-1} q|$ in 
(\ref{E:5.14}). The same argument with ${\rm div} \, \vec b$ in place of $q$ gives 
$$\int_{\R^n}  
  ( | \nabla (1-\Delta)^{-1} {\rm div} \, \vec b |^2 + |(1-\Delta)^{-1} \vec b|^2 ) \, |u|^2 \, dx 
 \le C \, ||u||^2_{W^{1, \, 2}(\R^n)}, 
$$
for all $u \in C^\infty_0(\R^n)$. The proof of Theorem~\ref{Theorem 5.1} is complete. 
\end{proof}

\section{Infinitesimal form boundedness and relative compactness}\label{Section 6} 

In this section, we discuss 
infinitesimal form boundedness and relative compactness properties (see \cite{RS}, \cite{Sch}) 
for the general second order differential operator $\mathcal{L}$. Since the coefficients of  $\mathcal{L}$ 
are arbitrary real- or complex-valued distributions, as above, we may assume without loss of generality 
that $\mathcal{L}$ is
in the divergence form $\mathcal{L} = -{\rm div} \, (A \nabla u) + \vec b \cdot \nabla +q$ where 
$A \in D'(\R^n)^{n\times n}$, $\vec b\in D'(\R^n)^n$, and $q \in D'(\R^n)$.   

The operator  $\mathcal{L}$ is said to be 
{\it relative form bounded\/}   with respect 
to the Laplacian on the (complex-valued)  
$L^2(\R^n)$ space  if 
\begin{equation}\label{E:6.1}
\left \vert \langle \mathcal{L} \, u, \, u  \rangle\right \vert 
\le \epsilon  \,  ||\nabla u||^2 _{L^{2} (\R^n)} + C(\epsilon) \, 
||u||^2_{L^{2} 
(\R^n)}, \quad 
u \in C^\infty_0(\R^n),
\end{equation} 
for some $\epsilon>0$ and $C(\epsilon)>0$.  This is obviously equivalent 
to the boundedness of the  sesquilinear form on 
$W^{1, \, 2} (\R^n)\times 
W^{1, \, 2} (\R^n)$, which was characterized in Theorem~\ref{Theorem 5.1}. 

However, in many applications it is of interest to distinguish the class 
of  $\mathcal{L}$ for which (\ref{E:6.1}) holds with relative bound zero, i.e., for 
every $\epsilon>0$ 
and some $C(\epsilon)>0$. In this case,   $\mathcal{L}$ 
 is said to be 
{\it infinitesimally form bounded\/} with respect to $-\Delta$ 
on $L^2(\R^n)$. 
For the potential energy operator $q \in D'(\R^n)$, the infinitesimal
 form boundedness with respect to $-\Delta$ was characterized recently in \cite{MV4}. Here we state the 
corresponding result for $\mathcal{L}$. 
 
Notice that, from the proof of Proposition \ref{Proposition 2.100} 
 applied to  (\ref{E:6.1}), it is immediate that the symmetric part  $A^s= \frac 1 2 (A+A^t)$ must 
 be equal to zero, while 
 the skew-symmetric part $A^c= \frac 1 2 (A-A^t)$ can be incorporated into $\vec b$ by letting 
 $\vec b_1 = \vec b - {\rm Div} \, A^c$, and considering {\rm (\ref{E:6.1})} for $\vec b_1 \cdot \nabla +q$. 
 Thus, without loss of generality it suffices to treat  the operator 
 $\mathcal{L} = \vec b \cdot \nabla +q$. 

\begin{theorem}\label{Theorem 6.1}
Let $\mathcal{L} = \vec b \cdot \nabla +q$, where $\vec b \in D'(\R^n)^n$ and $q \in D'(\R^n)$, $n \ge 2$. 
Then {\rm (\ref{E:6.1})} holds for every $\epsilon >0$ 
if and only if $\vec b$ and $q$ can be represented in 
the form {\rm (\ref{E:5.3})},  
where $F$ has vanishing mean oscillation, i.e., 
\begin{equation}\label{E:6.2} 
\lim_{\delta \to +0} \, 
\sup_{Q: \, |Q| \le \delta} \, \frac 1 {|Q|} \int_Q |F - m_Q(F)| \, dx =0,
\end{equation}
and 
\begin{equation}\label{E:6.3}
\lim_{\delta \to +0} \, \sup \, 
\left \{ \frac {\int_{Q}  \, 
|u|^2 \, d\mu} {||\nabla u||^2_{L^{2} (Q)}} : \quad 
u \in C^\infty_0(Q), \, \, u \not= 0, \quad |Q| \le \delta\right\} =0,
\end{equation}
where $d \mu = (|\vec c|^2 + |\vec h|^2 + |\gamma|) \, dx$. 
Moreover,  $\vec c $, $F$,  $\vec h$,  and $\gamma$ 
can be defined respectively by {\rm (\ref{E:5.5a})}, {\rm (\ref{E:5.5b})}.
\end{theorem}

The proof of Theorem~\ref{Theorem 6.1} follows by combining the approach of \cite{MV4}, which is based on a 
localization argument, with the form boundedness criterion obtained above.

\begin{remark} Analytic criteria for (\ref {E:6.3})
 to hold are discussed in \cite{MV4}.
 \end{remark} 

\begin{remark} Trudinger's condition where $C(\epsilon) = C \, \epsilon^{-\beta}$, 
$\beta >0$, 
in  (\ref {E:6.1}), and inequalities of Nash's type,
 \begin{equation}\label{E:6.1a}
\left \vert \langle \mathcal{L} \, u, \, u  \rangle\right \vert 
\le C  \,  ||\nabla u||^{2 \gamma} _{L^{2} (\R^n)}  \, 
||u||^{2(1-\gamma)}_{L^{1} 
(\R^n)}, \quad 
u \in C^\infty_0(\R^n),
\end{equation} 
where $\gamma\in (0, \, 1)$,  can be characterized using our approach as well; see \cite{MV4} where this is done for $\vec b =0$. 
 \end{remark}

Finally, we state a criterion  for the   {\it relative compactness} property which requires 
additional conditions at infinity. 

\begin{theorem}\label{Theorem 6.2}
Let $\mathcal{L} = \vec b \cdot \nabla +q$, where $\vec b \in D'(\R^n)^n$ and $q \in D'(\R^n)$, $n \ge 2$.
Then the operator $\vec b \cdot \nabla + q$ is relatively 
compact with respect to 
$-\Delta$ on $L^2(\R^n)$ 
if and only if $\vec b$ can be represented in 
the form {\rm (\ref{E:5.3})},  
where $F \in {\rm VMO}(\R^n)^{n \times n}$, i.e., 
\begin{align} 
& \lim_{\delta \to +0} \, \sup_{Q: \, |Q| \le \delta} \, 
\frac 1 {|Q|} \int_Q |F - m_Q(F)| \, dx =0,\\
& \lim_{\delta \to +\infty} \, \sup_{Q_0: \, |Q_0| \ge \delta} \, 
\frac 1 {|Q_0|} \int_{Q_0} |F - m_{Q_0} (F)| \, dx =0,
\end{align}
and 
\begin{align}
& \lim_{\delta \to +0} \, \sup \, 
\left \{ \frac {\int_{Q}  
|u|^2 \, d \mu} {|| u||^2_{W^{1, \, 2} (\R^n)}} : \quad 
u \in C^\infty_0(Q), \, \, u \not= 0, \quad |Q| \le \delta\right\} =0,\\
& \lim_{\delta \to +\infty} \, \sup \, 
\left \{ \frac {\int_{Q_0^c}  
|u|^2 \, d \mu} {|| u||^2_{W^{1, \, 2} (\R^n)}} : \quad 
u \in C^\infty_0(Q_0^c), \, \, u \not= 0, \quad |Q_0| \ge \delta\right\} =0,
\end{align}
where $Q_0$ denotes a cube  centered at the origin, and 
$d \mu = (|\vec c|^2 + |\vec h|^2 + |\gamma|) \, dx$. 
Moreover,  $\vec c$, $F$, $\vec h$,  and $q$ 
can be defined respectively by {\rm (\ref{E:5.5a})}, {\rm (\ref{E:5.5b})}.
\end{theorem}

The proof of Theorem~\ref{Theorem 6.2} is based on  the 
  form boundedness criterion obtained in the previous section, and is 
 analogous to the case  $\vec b =0$ treated in \cite{MV1}.

\end{document}